\documentclass{article}

\usepackage{amsmath}
\usepackage{amsthm}
\usepackage{graphicx}
\usepackage{amsfonts}
\usepackage{amssymb}
\usepackage{url}
\usepackage{amscd}
\usepackage{multicol}

\newcommand{\re}{\mathbb{R}}

\newcommand{\N}{\mathbb{N}}

\newcommand{\lmd}{\lambda}
\newcommand{\half}{\frac{1}{2}}

\newcommand{\dt}{\delta}
\newcommand{\Dt}{\Delta}
\newcommand{\af}{\alpha}
\newcommand{\bt}{\beta}

\newcommand{\sig}{\sigma}
\newcommand{\Sig}{\Sigma}

\newcommand{\supp}{\mbox{supp}}

\newcommand{\eproof}{{\ $\square$ }}
\newcommand{\reff}[1]{(\ref{#1})}

\newcommand{\prm}{\prime}

\newcommand{\bdes}{\begin{description}}
\newcommand{\edes}{\end{description}}

\newcommand{\bal}{\begin{align}}
\newcommand{\eal}{\end{align}}

\newcommand{\bnum}{\begin{enumerate}}
\newcommand{\enum}{\end{enumerate}}

\newcommand{\bit}{\begin{itemize}}
\newcommand{\eit}{\end{itemize}}

\newcommand{\bea}{\begin{eqnarray}}
\newcommand{\eea}{\end{eqnarray}}
\newcommand{\be}{\begin{equation}}
\newcommand{\ee}{\end{equation}}

\newcommand{\baray}{\begin{array}}
\newcommand{\earay}{\end{array}}

\newcommand{\bsry}{\begin{subarray}}
\newcommand{\esry}{\end{subarray}}

\newcommand{\bca}{\begin{cases}}
\newcommand{\eca}{\end{cases}}

\newcommand{\bcen}{\begin{center}}
\newcommand{\ecen}{\end{center}}

\newcommand{\bbm}{\begin{bmatrix}}
\newcommand{\ebm}{\end{bmatrix}}

\newcommand{\bmx}{\begin{matrix}}
\newcommand{\emx}{\end{matrix}}

\newcommand{\bpm}{\begin{pmatrix}}
\newcommand{\epm}{\end{pmatrix}}

\newcommand{\btab}{\begin{tabular}}
\newcommand{\etab}{\end{tabular}}

\newcommand{\thmlist}{
\begin{list}{Step 1}
{\setlength{\leftmargin}{0.6 in}\setlength{\labelwidth}{0.5 in}} }
\newcommand{\alglist}{
\begin{list}{Step 1}
{\setlength{\leftmargin}{1.1 in}\setlength{\labelwidth}{1.0 in}} }
\theoremstyle{plain}
\newtheorem{thm}{Theorem}[section]

\newtheorem{cor}[thm]{Corollary}

\newtheorem{ass}[thm]{Assumption}

\theoremstyle{definition}

\newtheorem{exm}[thm]{Example}
\newtheorem{alg}[thm]{Algorithm}

\newtheorem{remark}[thm]{Remark}

\renewcommand{\subsection}[1]{
    \stepcounter{subsection}
    \settowidth{\hangindent}{\bf\thesubsection.~}
    \hangafter=1
    \bigskip\bigskip\noindent
    {\bf\hbox{\thesubsection.~}#1}\par
    \nobreak
    \medskip
}

\usepackage[top=1.5in,bottom=1in,left=1.2in,right=1.2in]{geometry}

\begin{document}

\title{ Sparse SOS Relaxations for Minimizing Functions
that are Summations of Small Polynomials
\author{Jiawang Nie\footnote{
Department of Mathematics, University of California at San Diego,
9500 Gilman Drive, La Jolla, CA 92093.
\, Email:\,\mbox{njw@math.ucsd.edu} }
\,\, and\, James Demmel\footnote{Department of Mathematics and EECS, University of California,
Berkeley, CA 94720.\, Email:\,\mbox{demmel@cs.berkeley.edu} }
}
}

\maketitle

\abstract{ This paper discusses how to find the global minimum
of functions that are summations of small polynomials
(``small'' means involving a small number of variables).
Some sparse sum of squares (SOS) techniques are proposed.
We compare their computational complexity and lower bounds
with prior SOS relaxations.
Under certain conditions,
we also discuss how to extract the global minimizers
from these sparse relaxations.
The proposed methods are especially useful in
solving sparse polynomial system and nonlinear least squares problems.
Numerical experiments are presented,
which show that the proposed methods significantly improve
the computational performance of prior methods
for solving these problems.
Lastly, we present applications of this sparsity technique in
solving polynomial systems derived from nonlinear differential equations
and sensor network localization.
}

\bigskip

{\bf Key words:} Polynomials, sum of squares (SOS), sparsity, nonlinear least squares,
polynomial system, nonlinear differential equations,
sensor network localization

\section{Introduction}

Global optimization of multivariate polynomial functions
contains quite a broad class of optimization problems.
It has wide and important applications in science and engineering.
Recently, there has been much work on globally minimizing
polynomial functions using representation theorems
from real algebraic geometry
for positive polynomials.
The basic idea is to approximate nonnegative polynomials
by sum of squares (SOS) polynomials.
This approximation is also called {\it SOS relaxation}, since
not every nonnegative polynomial is SOS.
Here a polynomial is said to be SOS if
it can be written as a sum of squares of other polynomials.
The advantage of SOS polynomials is that
a polynomial is SOS if and only if
a certain semidefinite program (SDP)
formed by its coefficients is feasible.
Since SDP \cite{SDPhandbook} has efficient numerical methods,
we can check whether a polynomial is SOS
by solving a particular SDP.

To be more specific, suppose we wish to find the global
minimum value $f^*$ of a polynomial function $f(x)$ of
vector $x=(x_1,x_2,\cdots,x_n)\in \re^n$.
The SOS relaxation finds a lower bound $\gamma$ for $f^*$
such that the polynomial $f(x)-\gamma$ is SOS.
Obviously, $f(x)-\gamma$ being SOS implies that
$f(x)-\gamma$ is nonnegative for every real vector $x$.
Hence such a $\gamma$ is a lower bound.
The maximum $\gamma$ found this way
is called the SOS lower bound,
which is often denoted by $f^*_{sos}$.
The relation $f^*_{sos}\leq f^*$ always holds
for every polynomial $f(x)$
(it is possible that $f^*_{sos}=-\infty$).
When $f^*_{sos} = f^*$, we say the SOS relaxation
is {\it exact}.
We refer to \cite{Las01,PS01, ParMP} for more details
on SOS relaxations for polynomial optimization problems.
There are two important issues for applying SOS relaxation
in global optimization of polynomial functions:
the {\it quality} and {\it computational complexity}.

The quality means how good is the SOS lower bound $f^*_{sos}$.
In practice, as observed in \cite{PS01},
many nonnegative polynomials that are not
``artificially'' constructed are SOS.
However, Blekherman \cite{blek} pointed out that
there are much more nonnegative polynomials
than SOS polynomials.
But usually SOS relaxation provides very good
approximations, although theoretically it can fail
with high probability.
When SOS relaxation is not exact, i.e., $f^*_{sos} < f^*$,
there are methods to fix it by applying modified
SOS relaxations.
Nie, Demmel and Sturmfels \cite{njw_grad}
proposed to use SOS representations of $f(x)-\gamma$
modulo the gradient ideal of $f(x)$,
and show that the minimum value $f^*$
can be obtained when $f^*$ is attained at some point.
Schweighofer \cite{SchWg5} proposes to
minimize $f(x)$ over a semialgebraic set called
the {\it gradient tentacle},
and shows that the minimum value $f^*$
can be computed when $f^*>-\infty$ but not attainable.
Jibetean and Laurent \cite{JibeteanLaurent}, and
Lasserre \cite{Las0407} propose perturbing $f(x)$ by adding
a higher degree polynomial with tiny coefficients,
and showed that the lower bounds will
converge to the minimum value $f^*$.
Recently, Laurent \cite{LauSuv} gave a survey on solving polynomial optimization
by using semidefinite relaxations.
We refer to \cite{JibeteanLaurent,Las0407,LauSuv,njw_grad,SchWg5}
for related work.

Another important issue for SOS relaxation is the
computational complexity.
Suppose $f(x)$ has degree $2d$
(it must be even for $f(x)$ to have a finite minimum).
Then $f(x)$ has up to $\binom{n+2d}{2d}$ monomials.
The condition that $f(x)-\gamma$ being SOS
reduces to an SDP the size of whose linear matrix
inequality (LMI) is $\binom{n+d}{d}$ with
$\binom{n+2d}{2d}$ variables.
These numbers can be huge for
moderate $n$ and $d$, say, $n=2d=10$.
For large scale polynomial optimization problems,
the general SOS relaxation is very difficult to implement numerically.
Sometimes this complexity makes the applicability
of SOS relaxation very limited.
We refer to \cite{PS01,ParMP} for the connection between
SOS relaxation and SDP.

\medskip

\noindent
{\bf Prior work}\,
There is some work on exploiting sparsity in polynomial optimization
when the polynomials are {\it sparse}.
In such situations,
sparse SOS relaxations are available and
the resulting SDPs have reduced sizes, and hence larger
problems can be solved.
Here being sparse means that
the number of monomials with nonzero coefficients
is much smaller than the
maximum possible number $\binom{n+2d}{2d}$.
Kojima et al. \cite{Kjm03} and Parrilo \cite{ParStr03}
discussed how to exploit sparsity of SOS relaxations
in unconstrained polynomial optimization.
Kim et al. \cite{KKW} and Lasserre \cite{Las_sparse}
discussed sparse SOS relaxations
for constrained polynomial optimization problems
and showed convergence under certain conditions.
Waki et al. \cite{WKKM} proposed a heuristic procedure
to exploit sparsity for minimizing polynomials
by {\it chordal extension} of the
{\it correlation sparsity pattern graph} (csp graph):
the vertices of the csp graph are the variables
$x_1,\cdots,x_n$;
the edge $(x_i,x_j)$ exists whenever
$x_ix_j$ appears in one monomial of $f(x)$.
To find one chordal extension,
\cite{WKKM} proposed to use the symbolic sparse Cholesky
factorization of the csp matrix
with minimum degree ordering.
If the chordal extension of the csp graph is also sparse, 
then the sparsity technique in \cite{WKKM} works well.
However, if the chordal extension of the csp graph is much less sparse, 
then that sparsity technique might still
be too expensive to be implementable
for some practical problems.

%
%

\medskip
\noindent
{\bf Our contributions}\,
In many practical applications,
the polynomials are not only sparse,
but also given with certain sparsity patterns.
For instance, the polynomials
are often summations of other ``small'' polynomials,
i.e. polynomials only involving a small number of variables.
Sometimes, these representations contain
useful information that might help us save computations significantly.
These sparsity patterns are often ignored in prior work,
where these polynomials would be treated
using the usual ``dense'' algorithms.
The main contribution of this paper is to
propose new sparse SOS relaxation techniques
taking the given sparsity pattern into account,
and to show numerical experiments demonstrating
their accuracy and speed.

In this paper, we
consider the polynomial optimization problem of the form
\begin{align} \label{sumpoly}
\min_{x\in \re^n} \quad f(x)=\sum_{i=1}^m f_i(x_{\Dt_i})
\end{align}
where $\Dt_i \subset [n]=\{1,2,\cdots,n\}$.
Here each $f_i(x_{\Dt_i})$ is a
polynomial in $x_{\Dt_i}=(x_j|j\in \Dt_i)$.
Let $\deg(f_i)=2d_i$ and $2d=\deg(f)=\max\{2d_1,\cdots,2d_m\}$
(we assume each $f_i$ has even degree along with $f$).
One basic and natural idea for solving problem~\reff{sumpoly}
is to find the maximum $\gamma$ such that
\[ f(x) - \gamma = \sum_{i=1}^m s_i(x_{\Dt_i})\]
where each $s_i(x_{\Dt_i})$ is an SOS polynomial
in $x_{\Dt_i}$ instead of all the variables $x_1,\cdots,x_n$.
Exploiting this sparsity pattern can save significant computation
without sacrificing much solution quality
for many practical problems.
In addition to presenting its numerical implementation,
this paper will also
discuss the theoretical properties of this sparse SOS relaxation
and its variations.

The main distinction of our sparse SOS technique
from earlier work like Waki et al. \cite{WKKM}
is that we do not use the chordal extension of csp graphs.
In case that the csp graph of $f(x)$ in \reff{sumpoly}
is chordal, our sparsity technique is almost the same as
the one in \cite{WKKM}.
However, if the csp graph of $f(x)$ is not chordal
and its chordal extension is much less sparse,
then our sparsity technique is significantly more efficient.
If the csp graph of $f(x)$ is not chordal
and its chordal extension is also sparse,
then our sparsity technique is slightly more efficient
while not losing much quality of solution.
Furthermore, our sparsity technique can be applied to
solve bigger dense polynomial optimization problems
which can not be solved
by other existing methods.
This is due to the observation that every polynomial $g(x)$
is a summation of monomials whose number of variables
is at most the degree $\deg(g)$.
So, when $\deg(g)$ is small, like $4$ or $6$,
then the formulation~\reff{sumpoly} is a good sparse model.
The numerical computations show that
our sparsity technique is usually more efficient
than other existing methods
in solving problems of the form~\reff{sumpoly}.

Polynomial optimization problems of the form \reff{sumpoly}
have important practical applications:
(i) {\it Solving polynomial systems:}
Many large scale polynomial equations are often sparse
and each equation might involve just a few variables,
e.g., the polynomial equations obtained from discretization
in nonlinear differential equations.
Such polynomial systems can be equivalently transformed to
a global polynomial optimization problems of the form \reff{sumpoly}.
We will show that
the proposed sparse SOS relaxation is exact
when the polynomial system has at least one real solution.
(ii) {\it Nonlinear least squares:}
Many difficult problems in statistics, biology, engineering
or other applications require
solving certain nonlinear least squares problems
and finding their global optimal solutions.
If each equation is sparse, then
sparse polynomial optimization \reff{sumpoly}
is a very natural model
and our sparsity technique is very suitable.
Sensor network localization is one important application
of this kind.


\bigskip
\noindent
{\bf Outline}\,
This paper is organized as follows.
Section~2 introduces some notation and
background for SOS relaxations,
Section~3 discusses properties of the sparse SOS relaxation
and its variations,
Section~4 presents some numerical implementations,
and Section~5 shows applications.
Lastly, Section~6 draws some conclusions and discusses
future work in this area.

\section{ Preliminaries }
\setcounter{equation}{0}

This section introduces some notations and backgrounds
in SOS relaxation methods for minimizing polynomial functions.

Throughout this paper, we will use the following notation:
$\re$ is the field of real numbers;
$\N$ is the set of nonnegative integers;
$\re^{\Dt_i}$ =   $\{(x_{k_1},\cdots,x_{k_\ell}):\, x_{k_j}\in \re\}$
when $\Dt_i=\{k_1,k_2,\cdots,k_\ell\}$;
$\re[X]$: the ring of real polynomials in $X=(x_1,x_2,\cdots,x_n)$;
$\re[X_{\Dt_i}]$: the ring of real polynomials in $X_{\Dt_i}=(x_k)_{k\in\Dt_i}$;
$\sum\re[X]^2$: SOS polynomials in $\re[X]$;
$\sum\re[X_{\Dt_i}]^2$: SOS polynomials in $\re[X_{\Dt_i}]$;
$\sum\re_N[X]^2$: SOS polynomials in $\re[X]$ with degree at most $2N$;
$\sum\re_N[X_{\Dt_i}]^2$: SOS polynomials in $\re[X_{\Dt_i}]$ with degree at most $2N$;
$\|x\|_2$ =  $\sqrt{x_1^2+x_2^2+\cdots+x_n^2}$;
$x^\af$ =  $x_1^{\af_1}x_2^{\af_2}\cdots x_n^{\af_n}$ for $\af\in \N^n$;
$\supp(\af)$ =  $\{i\in [n]:\af_i \ne 0\}$;
$\supp(f)$ =  $\{\af\in \N^n: \text{ the coefficient of $x^\af$ in $f(x)$ is nonzero} \}$;
$|F|$ denotes the cardinality of set $F$;
$A^T$ denotes the transpose of matrix $A$;
$A\succeq (\succ) 0$ means matrix $A$ is positive semidefinite (definite);
$\mathcal{M}_d(y)$ is the moment matrix of order $d$ about $x\in \re^n$;
$\mathcal{M}_d^{\Dt_i}(y)$ is the moment matrix of order $d$ about $x_{\Dt_i}\in \re^{\Dt_i}$;
$\mathcal{M}_{F}(y)$ is the moment matrix generated monomials with support $F$.

\subsection{SOS and semidefinite programming (SDP) }

A polynomial $p(x)$ in $x=(x_1,\cdots,x_n)$
is said to be sum of squares (SOS) if
$p(x) =\sum_i p^2_i(x) $
for some polynomials $p_i(x)$.
Obviously, if $p(x)$ is SOS, then
$p(x)$ is {\it nonnegative}, i.e.,
$p(x)\geq 0$ for all $x\in \re^n$.
However, the converse is not true.
If $p(x)$ is nonnegative, then $p(x)$
is not necessarily SOS.
In other words, the set of SOS polynomials
(which forms a cone)
is properly contained in the set of
nonnegative polynomials
(which forms a larger cone).
The process of approximating nonnegative polynomials
by SOS polynomials is called SOS relaxation.
For instance, the polynomial
\begin{align*}
 &\, x_1^4+x_2^4+x_3^4+x_4^4-4x_1x_2x_3x_4 \\
= & \frac{1}{3}\left\{(x_1^2-x_2^2-x_4^2+x_3^2)^2+
(x_1^2+x_2^2-x_4^2-x_3^2)^2 +(x_1^2-x_2^2-x_3^2+x_4^2)^2+\right. \\
& \qquad \left. 2(x_1x_4-x_2x_3)^2+2(x_1x_2-x_3x_4)^2+2(x_1x_3-x_2x_4)^2 \right\}
\end{align*}
is SOS. This identity immediately implies that
\[ x_1^4+x_2^4+x_3^4+x_4^4-4x_1x_2x_3x_4 \geq 0,\,\,\,
\forall (x_1,x_2,x_3,x_4) \in \re^4,
\]
which is one arithmetic-geometric mean inequality.

The advantage of SOS polynomials over nonnegative polynomials
is that it is more tractable to check whether a polynomial is
SOS. To test whether a polynomial is SOS is equivalent to
testing the feasibility of some SDP \cite{PS01,ParMP},
which has efficient numerical solvers.
To illustrate this, suppose polynomial $p(x)$ has degree
$2d$ (SOS polynomials must have even degree).
Then $p(x)$ is SOS if and only if \cite{PS01,ParMP}
there exists a symmetric matrix $W\succeq 0$
such that
\[
p(x) = \mathbf{m}_d(x)^T\, W\, \mathbf{m}_d(x)
\]
where $\mathbf{m}_d(x)$ is the column vector of monomials up to
degree $d$. For instance,
\[
\mathbf{m}_2(x_1,x_2)=[\,1,\, x_1,\, x_2,\, x_1^2,\,  x_1x_2,\,  x_2^2\,]^T.
\]
As is well-known, the number of monomials in $x$
up to degree $d$ is
$\binom{n+d}{d}$. Thus the size of matrix $W$
is $\binom{n+d}{d}$.
This number can be very large.
For instance, when $n=d=10$,
$\binom{n+d}{d}\geq 10^{5}$.
However, for fixed $d$ (e.g.,$d=2$),
$\binom{n+d}{d}$ is polynomial in $n$.
On the other hand, it is NP-hard (with respect to $n$) to tell whether
a polynomial is nonnegative
whenever $2d\geq 4$ (even when $d$ is fixed)\cite{Las01}.

\subsection{ SOS relaxation in polynomial optimization}

Let $f(x)=\sum_{\alpha }\, f_\alpha x^\alpha$
be a polynomial in $x$.
Consider the global optimization problem
\[
f^*:=\min_{x\in \re^n}\,\ f(x).
\]
This problem is NP-hard when $\deg(f)\geq 4$.
The standard SOS relaxation is
\begin{align*}
f^*_{sos}:=\max & \quad \gamma \\
s.t. & \quad f(x)-\gamma \text{ is SOS }.
\end{align*}
Obviously we have that $f^*_{sos}\leq f^*$.
In practice, SOS provides very good approximations,
and often gives exact global minimum,
i.e., $f^*_{sos}=f^*$,
even though theoretically there are many more
nonnegative polynomials than SOS polynomials \cite{blek}.

In terms of SDP, the SOS relaxation can also be written as
\begin{align}
f^*_{sos}:=\max & \quad \gamma  \label{sos1}\\
s.t. & \quad f(x)-\gamma =
\mathbf{m}_d(x)^T W \mathbf{m}_d(x) \label{sos2} \\
&\quad W \succeq 0 \label{sos3}
\end{align}
where $2d=\deg(f)$.
The decision variable in the above
is $(\gamma,W)$ instead of $x$.
The above program is convex about $(\gamma, W)$.
A lower bound $f^*_{sos}$ can be computed
by solving the resulting SDP.
It can be shown \cite{Las01} that the dual of \reff{sos1}-\reff{sos3} is
\begin{align}
f^*_{mom}:=
\min_{y} & \quad \sum_{|\alpha|\leq 2d} f_\alpha y_\alpha  \label{mom1}\\
s.t. & \quad \mathcal{M}_d(y) \succeq 0 \label{mom2} \\
& \quad y_{0,\cdots,0}=1. \label{mom3}
\end{align}
Here $\mathcal{M}_d(y)$ is the {\it moment matrix} generated by
$y=(y_\alpha)$, a vector indexed by monomials of degree at most $2d$.
The rows and columns of moment matrix $\mathcal{M}_d(y)$
are indexed by integer vectors.
Each entry of $\mathcal{M}_d(y)$ is defined as
\[
\mathcal{M}_d(y)(\alpha,\beta):=y_{\alpha+\beta},\,\, \forall
|\alpha|,|\beta| \leq d.
\]
For instance, when $d=2$ and $n=2$, the vector
\[
y =[\, y_{0,0},\, y_{1,0},\, y_{0,1},\, y_{2,0},\, y_{1,1},\, y_{0,2},\,
y_{3,0},\, y_{2,1},\, y_{1,2},\, y_{0,3},\,
y_{4,0},\, y_{3,1},\, y_{2,2},\, y_{1,3},\, y_{0,4}\,],
\]
defines moment matrix
\[
\mathcal{M}_2(y) = \bbm
y_{0,0} & y_{1,0} & y_{0,1} & y_{2,0} & y_{1,1} & y_{0,2} \\
y_{1,0} & y_{2,0} & y_{1,1} & y_{3,0} & y_{2,1} & y_{1,2} \\
y_{0,1} & y_{1,1} & y_{0,2} & y_{2,1} & y_{1,2} & y_{0,3} \\
y_{2,0} & y_{3,0} & y_{2,1} & y_{4,0} & y_{3,1} & y_{2,2} \\
y_{1,1} & y_{2,1} & y_{1,2} & y_{3,1} & y_{2,2} & y_{1,3} \\
y_{0,2} & y_{1,2} & y_{0,3} & y_{2,2} & y_{1,3} & y_{0,4}
\ebm.
\]
For SOS relaxation \reff{sos1}-\reff{sos3} and
its dual problem \reff{mom1}-\reff{mom3},
strong duality holds \cite{Las01}, i.e.,
their optimal values are equal ($f^*_{sos}=f^*_{mom}$).
Hence $f^*_{mom}$ is also a lower bound for
the global minimum $f^*$ of $f(x)$.

Now let us see how to extract minimizer(s)
from optimal solutions to \reff{mom1}-\reff{mom3}.
Let $y^*$ be one optimal solution.
If moment matrix $\mathcal{M}_d(y^*)$ has rank one, then
there exists one vector $w$ such that
$\mathcal{M}_d(y^*)=ww^T$. Normalize $w$ so that $w_{(0,\cdots,0)}=1$.
Set $x^*=w(2:n+1)$. Then the relation $\mathcal{M}_d(y^*)=ww^T$
immediately implies that
$y^*=\mathbf{m}_{2d}(x^*)$, i.e., $y^*_\alpha = (x^*)^\alpha$,
so $f^*_{mom}=f(x^*)$.
This says that a lower bound of $f(x)$
is attained at one point $x^*$.
So $x^*$ is one global minimizer.

When moment matrix $\mathcal{M}_d(y^*)$ has rank
more than one, the process described above
does not work. However, if $\mathcal{M}_d(y^*)$ satisfies
the so-called {\it flat extension condition }
\[
\mbox{rank}\, \mathcal{M}_{k}(y^*) = \mbox{rank}\, \mathcal{M}_{k+1}(y^*)
\]
for some $0\leq k \leq m-1$,
we can extract more than one minimizer
(in this case the global solution is not unique).
When the flat extension condition is met, it can be shown \cite{CurFia}
that there exist distinct vectors $u_1,\cdots,u_r$ such that
\[
\mathcal{M}_d(y^*)=\lmd_1 \mathbf{m}_d(u_1)\cdot \mathbf{m}_d(u_1)^T + \cdots +
\lmd_r \mathbf{m}_d(u_r)\cdot \mathbf{m}_d(u_r)^T
\]
for some $\lmd_i>0,\, \sum_{i=1}^r\lmd_i=1$.
Here $r=\mbox{rank}\, \mathcal{M}_d(y^*)$.
The set $\{u_1,\cdots,u_r\}$ is called an {\it $r$-atomic representing support}
for moment matrix $\mathcal{M}_d(y^*)$.
All the vectors $u_1,\cdots,u_r$ can be shown to
be global minimizers.
They can be computed by solving some
particular eigenvalue problem.
We refer to \cite{CurFia} for flat extension conditions
in moment problems and \cite{HenLas}
for extracting minimizers.

\subsection{Exploiting sparsity in SOS relaxation}

As mentioned in the previous subsections,
the size of matrix $W$ in SOS relaxation
is $\binom{n+d}{d}$,
which can be very large.
So SOS relaxation is expensive
when either $n$ or $d$ is large.
This is true for general dense polynomials.
However, if $f(x)$ is sparse, i.e.,
its support $\mathcal{F}=\mbox{supp}(f)$ is small,
the size of the resulting SDP can be reduced significantly.
Without loss of generality, assume $(0,\cdots,0)\in \mathcal{F}$.
Then $\mbox{supp}(f)=\mbox{supp}(f-\gamma)$ for any number $\gamma$.

\bigskip

Suppose $f(x)-\gamma =\sum_i \phi_i(x)^2$ is an SOS decomposition.
Then by Theorem~1 in \cite{Rez78} we have
\[
\mbox{supp}(\phi_i)\subset \mathcal{F}^0 :=
 \left( \text{ the convex hull of }
\half \mathcal{F}^e \right)
\]
where $\mathcal{F}^e=\{\alpha \in \mathcal{F}:\,\,
\alpha \, \text{ is an even integer vector} \}$.
There exist some work \cite{Kjm03,WKKM}
on exploiting sparsity further.
Here we briefly describe the technique introduced in \cite{WKKM}.

For polynomial $f(x)$, define its {\it csp graph}
$G=([n],E)$ such that $(i,j)\in E$ if and only if
$x_ix_j$ appears in some monomial of $f(x)$.
Let $\{C_1,C_2,\cdots,C_K\}$ be the set of all
maximal cliques of graph $G$.
Waki et al. \cite{WKKM} proposed to represent $f(x)-\gamma$ as
\[
f(x)-\gamma = \sum_{i=1}^K s_i(x),\quad
\text{ each } s_i(x) \text{ being SOS }
\mbox{supp}(s_i) \subset C_i.
\]
Theoretically, when $f(x)-\gamma$ is SOS,
the above representation may not hold
(see Example~\ref{ripnotsuf}).
And it is also difficult to find all the maximal cliques of
graph $G$.
Waki et al. \cite{WKKM} propose to replace
$\{C_1,C_2,\cdots,C_K\}$
by the set of all maximal cliques of one {\it chordal extension} of $G$.
We refer to \cite{BlaPey93} for properties of chordal graphs.
For chordal graphs, there are efficient methods to
find all the maximal cliques.
Chordal extension is essentially
the {\it sparse symbolic Cholesky factorization}.
See \cite{WKKM} for more details on how to get the
chordal extension.

We remark that in the worst case the sparse SOS relaxation above
might be weaker than the general dense SOS relaxation
even when the chordal extension is applied,
as shown by Example~\ref{ripnotsuf}.

There is much work in
exploiting sparsity in SOS relaxations.
We refer to \cite{GatPar,KKW,Kjm03,ParStr03,WKKM}
and the references therein.

\section{ The sparse SOS relaxation }
\setcounter{equation}{0}

Throughout this paper, we assume
$
f(x) = \sum_{i=1}^m \, f_i(x_{\Dt_i}).
$
Let $\|\Dt\|$ be the maximum cardinality of $\Dt_i$, i.e.,
$
\|\Dt\| = \max_i \, |\Dt_i|.
$
We are interested in the case that $\|\Dt\| \ll n$.
To find the global minimum $f^*$ of $f(x)$,
we propose the following sparse SOS relaxation
\begin{align*}
f^*_\Dt :=\max & \quad \gamma \\
s.t. & \quad f(x) - \gamma \in
\sum_{i=1}^m \sum \re_{d}[x_{\Dt_i}]^2.
\end{align*}
In terms of SDP, the above SOS relaxation is
essentially the same as
\begin{align}
f^*_\Dt :=\max & \quad \gamma \label{sumsos1}\\
s.t. & \quad f(x) - \gamma = \sum_{i=1}^m
\mathbf{m}_d(x_{\Dt_i})^T W_i \mathbf{m}_d(x_{\Dt_i})  \label{sumsos2} \\
& \qquad W_i\succeq 0, \, i=1,\cdots,m. \label{sumsos3}
\end{align}
Notice that \reff{sumsos2} is an identity.
Let
\be \label{def:F}
\mathcal{F}_i=\{ \af\in \N^n:
\supp(\af) \subset \Dt_i,\, |\af| \leq 2d\}, \,\,\,\,
\mathcal{F}=\bigcup\mathcal{F}_i.
\ee
Write $f(x)=\sum_\af f_\af x^\af$.
Since $f(x)=\sum_i f_i(x_{\Dt_i})$,
$f_\af\ne 0$ implies that $\af\in \mathcal{F}$.
By comparing coefficients of both sides of \reff{sumsos2},
we have equality constraints
\begin{align} \label{coefeq}
f_0 -\gamma = \sum_{i=1}^m \,  W_i(0,0), \quad
f_\af = \sum_{i=1}^m \, \sum_{\eta+\tau =\af} W_i(\eta,\tau), \,\,\,
\forall \, \af \ne 0.
\end{align}

Now we derive the dual problem for
\reff{sumsos1}-\reff{sumsos3}.
Notice that constraint \reff{sumsos2} is equivalent to the equality constraints
\reff{coefeq}. Let $y=(y_\af)_{\af\in \mathcal{F}}$
be the Lagrange multipliers
for equations in \reff{coefeq}, and $U_i$ be the Lagrange multipliers
for inequalities in \reff{sumsos3}. Each $U_i$ is also positive semidefinite.
The Lagrange function for problem
\reff{sumsos1}-\reff{sumsos3} is
\begin{align*}
\mathcal{L}&= \gamma +
\left(f_0-\gamma-\sum_i W_i(0,0)\right)y_0
+\sum_{ 0 \ne \af \in \mathcal{F}}
\left(f_\af-\sum_{i=1}^m \sum_{\eta+\tau =\af} W_i(\eta,\tau) \right)y_\af +
\sum_i W_i\bullet U_i \\
&= \gamma(1-y_0)+\sum_{\af\in \mathcal{F}} f_\af y_\af +
\sum_{i=1}^m \sum_{\af\in \mathcal{F} } \sum_{\eta+\tau=\af}
W_i(\eta,\tau) \Big(U_i(\eta,\tau) - y_\af \Big).
\end{align*}
So we can see that
\[
\max_{\gamma, W_i} \,\mathcal{L}(\gamma,W_i,y_\af, U_i) =
\bca
\sum_\af f_\af y_\af  & \text{ if } y_0=1, U_i=\mathcal{M}_d^{\Dt_i}(y) \succeq 0 \, ;\\
+\infty & \text{ otherwise.}
\eca
\]
Therefore the dual of
\reff{sumsos1}-\reff{sumsos3} is
\begin{align}
f^*_\Sig:=
\min & \quad  \sum_{\af\in \mathcal{F}} f_\af y_\af \label{summom1}\\
s.t. & \quad  \mathcal{M}_d^{\Dt_i}(y) \succeq 0,
\, i=1,\cdots,m  \label{summom2} \\
& \quad y_{0} =1. \label{summom3}
\end{align}

\subsection{Complexity comparison}

Since the dual of the standard or sparse SOS relaxation
not only returns the SOS lower bound but also
provides the moment matrix to help extract minimizers,
we compare the computational complexity
of \reff{mom1}-\reff{mom3} and \reff{summom1}-\reff{summom3}.
The LMI~\reff{mom2} is of size $\binom{n+d}{d}=\mathcal{O}(n^d)$
and has $\binom{n+2d}{2d}=\mathcal{O}(n^{2d})$ decision variables.
At each step of an interior-point method (e.g., the dual scaling method \cite{DSDP3}),
the complexity for solving \reff{mom1}-\reff{mom3}
is $\mathcal{O}(n^{6d})$.
On the other hand, \reff{summom2} has $m$ LMIs, which are of sizes
at most $\binom{\|\Dt\|+d}{d}=\mathcal{O}(\|\Dt\|^d)$,
and $\mathcal{O}\left(m\binom{\|\Dt\|+2d}{2d}\right)
=\mathcal{O}(m\|\Dt\|^{2d})$ decision variables.
At each step of interior-point methods,
the complexity for solving \reff{summom1}-\reff{summom3}
is $\mathcal{O}(m^3\|\Dt\|^{6d})$.
When $\|\Dt\|$ is independent of $n$ and
$m=\mathcal{O}(n^p)$ with $p<2d$, then
\[
\mathcal{O}(m^3\|\Dt\|^{6d})\ll \mathcal{O}(n^{6d}).
\]
Therefore \reff{summom1}-\reff{summom3} is much easier to solve
than \reff{mom1}-\reff{mom3}.


The complexity of sparse SOS relaxation in \cite{WKKM} depends on
the chordal extension of the csp graph.
In the worst case, it can be as big as for the general
SOS relaxation \reff{mom1}-\reff{mom3}.
Let $\Omega$ be the maximum size of the maximal cliques of the chordal extension.
In practice, $\Omega$ is often bigger than or equal to $\|\Dt\|$.
When $\Omega > \|\Dt\|$, the SOS relaxation \reff{summom1}-\reff{summom3}
is usually more efficient.

\subsection{Lower bound analysis}

Recall that $\mathcal{F}_i=\{ \af\in \N^n:
\supp(\af) \subset \Dt_i,\, |\af| \leq 2d\}$. From the
representation~\reff{sumpoly} of $f(x)$, we have
\[
\supp(f) \subseteq \bigcup_{i=1}^m\, \mathcal{F}_i.
\]
This leads us to think that
the relaxation \reff{summom1}-\reff{summom3}
should give reasonable lower bounds,
although it might be weaker than the general SOS
(see Example~\ref{ripnotsuf}).

\begin{thm} \label{bdrel}
The optimal values $f^*_{\Sig},f^*_\Dt,f^*_{sos}, f^*$
satisfy the relationship
\[
f^*_\Sig=f^*_\Dt \leq f^*_{sos} \leq f^*.
\]
\end{thm}
\proof
The latter two inequalities are obvious
because the feasible region defined by
\reff{summom2}-\reff{summom3}
contains the one defined by
\reff{mom2}-\reff{mom3}.
To prove the first equality,
by the standard duality argument for convex program,
it suffices to show that
\reff{summom2} admits a strict interior point.
Define $\hat y = (\hat y_\af)_{\af \in \mathcal{F}}$ as
\[
\hat y_\alpha := \frac{\int_{\re^{n}} x^\alpha
e^{ -\|x\|_2^2}\, d x}{ \int_{\re^{n}} e^{-\|x\|_2^2}\, dx}.
\]
For every nonzero vector
$\xi=(\xi_\alpha)_{\af \in \mathcal{F}_i}$, we have
\[
\xi^T\,M_d^{\Dt_i}(\hat y)\,\xi = \frac{
\int_{\re^n}
\left(\sum_{|\alpha|\leq d} \xi_\af x^\af \right)^2
e^{-\|x\|_2^2} d x}
{\int_{\re^n} e^{ -\|x\|_2^2} d x } >0.
\]
So $M_d^{\Dt_i}(\hat y)\succ 0$ for every $1\leq i \leq m$.
Therefore $\hat y$ is an interior point for
\reff{summom1}-\reff{summom3},
which implies the strong duality $f^*_\Sig=f^*_\Dt$.
\eproof
\begin{remark}
Theorem~\ref{bdrel} implies that the lower bound
$f^*_\Dt$ given by \reff{sumsos1}-\reff{sumsos3}
is weaker than the SOS lower bound $f^*_{sos}$.
There are examples such that
$f^*_\Dt<f^*_{sos}$ (see Example~\ref{ripnotsuf}).
However, in many numerical simulations,
the lower bound $f^*_\Dt$ is very useful.
For randomly generated polynomials,
as shown in Section~4,
it frequently happens that $f^*_\Dt=f^*_{sos}$.
On the other hand, under some conditions,
we can prove $f^*_\Dt=f^*_{sos}$.
\end{remark}
\bigskip

Suppose $\Dt_1,\Dt_2,\cdots,\Dt_m$ satisfy the
{\it running intersection property}:
\begin{align}\label{rip}
\mbox{ For every }\, 1\leq i\leq m-1,\,\,
\exists\, k\leq i \,\, \mbox{ such that } \quad
\Dt_{i+1} \cap \left(\bigcup_{j=1}^i\Dt_j \right)  \subsetneq \Dt_k.
\end{align}

\begin{thm}  \label{thm:Dt=sos}
Suppose \reff{summom1}-\reff{summom3} has a
optimal solution $y^*$ such that
each $\mathcal{M}_{d_i}^{\Dt_i}(y^*)$ has
a representing measure $\mu_i$ on $\re^{\Dt_i}$.
If condition~\reff{rip} holds,
then $f^*_\Dt=f^*_{sos}$.
\end{thm}
\proof
For any $\Dt_i,\Dt_j$,
$\mathcal{M}_d^{\Dt_i\cap \Dt_j}(y^*)$ is
a common principle submatrix of $\mathcal{M}_d^{\Dt_i}(y^*)$
and $\mathcal{M}_d^{\Dt_j}(y^*)$.
So the marginals of measures $\mu_i$
are consistent, i.e.,
the restrictions of these measures
on the common subspaces are the same.
By Lemma~6.4 in \cite{Las_sparse},
there exists a measure on $\re^n$ such that
$\mu_i$ is the marginal of $\mu$
with respect to $\Dt_i$ for all $i=1,\cdots,m$.
Define vector $\tilde y$ such that
\[
\mathcal{M}_d(\tilde y)=\int_{\re^n} \,
\mathbf{m}_d(x) \mathbf{m}_d(x)^T \mu(dx).
\]
Then every $\mathcal{M}_d^{\Dt_i}(y^*)$
is a principle submatrix of $\mathcal{M}_d(\tilde y)$.
So $\tilde y_\af=y^*_\af$ whenever
$\supp(\af)\subset \Dt_j$ for some $j$.
Since the $f_\af \ne 0$ implies
$\supp(\af)\subset \Dt_j$ for some $j$,
we know the objective value of \reff{summom1}
is the same for $y^*$ and $\tilde y$.
Thus $f^*_{sos}\leq f^*_\Dt$.
Since $f^*_{sos}\geq f^*_\Dt$,
we get $f^*_{sos} = f^*_\Dt$.
\eproof

\begin{remark}
The running intersection property~\reff{rip} alone
is not sufficient to guarantee the equality $f^*_\Dt=f^*_{sos}$,
as shown by the following example.
\end{remark}

\begin{exm} \label{ripnotsuf}
$f(x)=f_1(x_1,x_2)+f_2(x_2,x_3)$
where $f_1=x_1^4+(x_1x_2-1)^2$ and $f_2=x_2^2x_3^2+(x_3^2-1)^2 $.
Solving dense SOS relaxation~\reff{sos1}-\reff{sos3}
and sparse SOS relaxation~\reff{sumsos1}-\reff{sumsos3}
numerically, we find that
\[
f^*_\Dt \approx 5.0\cdot 10^{-5} \, < \,f^*_{sos}\approx 0.8499.
\]
Actually the minimum $f^*\approx 0.8650$.
First, solve equation $\nabla f(x)=0$,
and evaluate $f(x)$ on these critical points,
then we find the minimum of these critical values is about $0.8650$.
So $f^*<1$.
Second, we prove that the minimum $f^*$ is attainable.
Let $\{x^{(k)}\}$ be a sequence such that
$f(x^{(k)})\to f^*$ as $k$ goes to infinity.
We claim that the sequence $\{x^{(k)}\}$ must be bounded.

\begin{minipage}[c]{0.75\textwidth}
Otherwise, suppose $x^{(k)}\to \infty$. Thus at least one
of coordinates $x^{(k)}_1,\,x^{(k)}_2,\,x^{(k)}_3$ should
go to infinity. If either $x^{(k)}_1$ or $x^{(k)}_3$ goes to infinity,
then $f(x^{(k)})$ goes to infinity, which is not possible.
So $x^{(k)}_2\to \infty$.
Since $\{f(x^{(k)})\}$ is bounded,
without loss of generality,
we assume $x^{(k)}_1\to a_1,\,x^{(k)}_1x^{(k)}_2\to a_{12},\,
x^{(k)}_2x^{(k)}_3\to a_{23},\,  x^{(k)}_3\to a_3$
for some numbers $a_1,a_{12},a_{23},a_3$.
If $a_3 = 1$, then $x^{(k)}_2$ is convergent to $a_2$,
which is not possible. And, if $a_3\ne 0$, then
$x^{(k)}_2x^{(k)}_3$ goes to infinity, which is also not possible.
So $a_3=0$, and hence
\[
f(x^{(k)}) \geq (x_3^2-1)^2 \to 1 > f^*,
\]
which is a contradiction.
\end{minipage}

\noindent
So the sequence $\{x^{(k)}\}$ is bounded and has an accumulation point
$x^*$. Then we must have $f(x^*)=f^*$, which means that $f^*$ is attained
at some point. From the computation of critical values, we know
$f^*\approx 0.8650$. For this polynomial,
both the dense and sparse SOS relaxation are not exact:
$ f^*_\Dt < f^*_{sos} <f^*$,
and the method in \cite{WKKM} gives the same lower bound $f^*_\Dt$.
\end{exm}

\begin{cor}~\label{quadpoly}
If all $f_i$ are quadratic and condition~\reff{rip}
holds, then
$f^*_{sos}=f^*_\Dt$.
\end{cor}
\proof
When all $f_i$ are quadratic, i.e., $d_i=1$,
the entries of moment matrix $\mathcal{M}^{\Dt_i}_1$
are the first and second order moments.
The positive semidefiniteness of $\mathcal{M}^{\Dt_i}_1$
implies $\mathcal{M}^{\Dt_i}_1$
has a representing measure.
Then the conclusion is immediately
implied by Theorem~\ref{thm:Dt=sos}.
\eproof

\begin{remark}
If the running intersection condition~\reff{rip} fails,
then Corollary~\ref{quadpoly} is no longer true,
as shown by the example below.
\end{remark}

\begin{exm} Consider the polynomial
$f(x)=f_1(x_1,x_2)+f_2(x_2,x_3)+f_3(x_1,x_3)$
where $f_1=\half(x_1^2+x_2^2)+2x_1x_2$, $f_2=\half(x_2^2+x_3^2)+2x_2x_3$,
$f_3=\half(x_1^2+x_3^2)+2x_1x_3$. In this case
\[
\Dt_1=\{1,2\},\quad  \Dt_2=\{2,3\},\quad  \Dt_3=\{1,3\}.
\]
The running intersection property~\reff{rip} fails. But we have
$f^*_\Dt=-\infty < f^*_{sos}=f^*=0$.
\end{exm}

\subsection{ Extraction of minimizers}

In this subsection, we discuss
how to extract minimizer(s) $x^*=(x_1^*,\,\cdots,\,x_n^*)$.
Suppose $y^*=(y^*_\af)_{\af\in\mathcal{F}}$ is one optimal solution to
\reff{summom1}-\reff{summom3}. Let $\dt_i=\{i\}$ for every $i$.
The entries of $y^*$ whose indices are supported in $\dt_i$ are
\[
y^*_{0},\, y^*_{e_i},\, y^*_{2e_i},\, \cdots,\, y^*_{2de_i}
\]
which are the entries the moment matrix $M_d^{\dt_i}(y^*)$.
So coordinate $x_i^*$ can be extracted
from moment matrix $M_d^{\dt_i}(y^*)$
if it satisfies the flat extension condition.
Let $\mathcal{V}_i$ be the set of all the points
that can be extracted from the moment matrix
$M_d^{\dt_i}(y^*)$.
If $\mathcal{V}_i$ is a singleton,
then $x_i^*$ has a unique choice.

The situation is more subtle if some $\mathcal{V}_i$ has
cardinality greater than one.
Suppose for some $i,j\in [n]$ we have
$|\mathcal{V}_i|>1$ and $|\mathcal{V}_j|>1$.
Can $x_i^*\,x_j^*$ appear simultaneously
in the optimal solution $x^*$ for arbitrarily chosen
$x_i^*\in\mathcal{V}_i,\, x_j^*\in \mathcal{V}_j$?
The answer is obviously no!
For instance, the polynomial
\[ (x_1^2-1)^2+(x_2^2-1)^2+(x_1-x_2)^2\]
has only two global minimizers $\pm (1,1)$.
We find that
$
\mathcal{V}_1=\mathcal{V}_2=\left\{ 1,-1 \right\}.
$
But obviously $(1,-1)$ and $(-1,1)$ are not global minimizers.

Now what is the rule for matching $x_i^*$ and $x_j^*$
if  $|\mathcal{V}_i|>1$ or $|\mathcal{V}_j|>1$?
So far we have not yet used the information of moment matrix
$M_d^{\Dt_i}(y^*)$. If $M_d^{\Dt_i}(y^*)$
also satisfies the flat extension condition,
we can extract the tuples $x^*_{\Dt_i}=(x_k^*)_{k\in\Dt_i}$
from $M_d^{\Dt_i}(y^*)$.
Let $\mathcal{X}_{\Dt_i}$ be set of all such tuples
that can be extracted from $M_d^{\Dt_i}(y^*)$.
One might ask whether $\mathcal{V}_i$ and
$\mathcal{X}_{\Dt_i}$ are consistent, that is,
does $x^*_{\Dt_i} \in \mathcal{X}_{\Dt_i}$ imply that
$x_k^*\in \mathcal{V}_k$ for all $k\in \Dt_i$?
Under the flat extension assumption,
the answer is yes, which
is due to the following theorem.

\begin{thm}
Suppose $y^*$ is one optimal solution to
\reff{summom1}-\reff{summom3} such that
all $M_d^{\Dt_i}(y^*)$
satisfy the flat extension condition.
Then for any
$x^*_{\Dt_i} \in \mathcal{X}_{\Dt_i}$, it holds that
$x_k^*\in \mathcal{V}_k$ for all $k\in \Dt_i$.
\end{thm}
\proof
Let $\mathcal{X}_{\Dt_i}= \{
x^{(1)}_{\Dt_i},\, x^{(2)}_{\Dt_i},\,\cdots,\, x^{(r)}_{\Dt_i}
\}$
be the $r$-atomic representing support for $\mathcal{M}_d^{\Dt_i}(y^*)$.
Then we have decomposition
\[
M_d^{\Dt_i}(y^*) = \sum_{\ell=1}^r \,\lmd_\ell\,
\mathbf{m}_2(x^{(\ell)}_{\Dt_i})
\mathbf{m}_2(x^{(\ell)}_{\Dt_i})^T
\]
for some $\lmd_1,\cdots,\lmd_r>0,\sum_{\ell=1}^r \,\lmd_\ell=1$.
Notice that $M_d^{\dt_k}(y^*)$ is a principle submatrix of
$M_d^{\Dt_i}(y^*)$. So we also have that for every $k\in \Dt_i$
\[
M_d^{\dt_k}(y^*) = \sum_{\ell=1}^r \,\lmd_\ell\,
\mathbf{m}_2( x_k^{(\ell)} )
\mathbf{m}_2( x_k^{(\ell)} )^T.
\]
This means that
$\{x^{(1)}_k,\, x^{(2)}_k,\,\cdots,\, x^{(r)}_k \}$
is a $r$-atomic representing support for moment matrix
$M_d^{\dt_k}(y^*)$ (some $x_k^{(\ell)}$ might be the same).
By the definition of $\mathcal{V}_i$, we have
$\{x_k^{(1)},\cdots,x_k^{(r)}\}\subseteq \mathcal{V}_k$.
\eproof

\begin{thm} \label{optmin}
Suppose $y^*$ is one optimal solution to
\reff{summom1}-\reff{summom3} such that
all $M_d^{\Dt_i}(y^*)$ satisfy the flat extension condition.
Then any $x^*=(x_1^*,\,\cdots,\,x_n^*)$
with $x_k^*\in \mathcal{V}_k$ and $x^*_{\Dt_i}\in \mathcal{X}_{\Dt_i}$
for all $k$ and $i$
is a global optimal minimizer of $f(x)$.
\end{thm}
\proof
Fix $x^*$ as in the theorem.
Since $M_d^{\Dt_i}(y^*)$ satisfies the flat extension condition,
we have the decomposition
\[
M_d^{\Dt_i}(y^*) =
\lmd_{\Dt_i} \mathbf{m}_d(x^*_{\Dt_i} )
\mathbf{m}_2( x^*_{\Dt_i})^T+\hat M_{\Dt_i}
\]
where $1\geq \lmd_{\Dt_i} >0$ and
$\hat M_{\Dt_i}\succeq 0$.
Now let
$\lmd=\min_{i}\lmd_{\Dt_i}>0$
and
\[
M_{\Dt_i} =
(\lmd_{\Dt_i}-\lmd) \mathbf{m}_2(x^*_{\Dt_i})
\mathbf{m}_2(x^*_{\Dt_i})^T
+\hat M_{\Dt_i} \succeq 0.
\]
Notice that $\hat M_{\Dt_i}$ and $M_{\Dt_i} $ are also moment matrices.
Without loss of generality, we can assume
$\lmd <1$, since otherwise each $M_2^{\Dt_i}(y^*)$
has rank one and then $x^*$ is obviously a global minimizer.
For every $\af \in \mathcal{F}_i$, define
$ \hat y_{\af}= \left(x_{\Dt_i}^*\right)^{\af}$
and $\hat y = (\hat y_\af)_{\af\in \mathcal{F}}$.
Let $\tilde y =(\tilde y_\af)_{\af\in \mathcal{F}}$
be the vector such that $y^* = \lmd \hat y +(1-\lmd)\tilde y$.
Then it holds
\[
M_d^{\Dt_i}(y^*) = \lmd  M_d^{\Dt_i}(\hat y)+(1-\lmd)M_d^{\Dt_i}(\tilde y).
\]
Obviously vector $\tilde y$ is feasible for
\reff{summom2}-\reff{summom3} since
\[
M_d^{\Dt_i}(\tilde y)= \frac{1}{1-\lmd}
\left( M_d^{\Dt_i}(y^*) - \lmd  M_d^{\Dt_i}(\hat y) \right)
= \frac{1}{1-\lmd} M_{\Dt_i} \succeq 0.
\]
Since $y^*$ is optimal, we can see
$\sum_{\af\in \mathcal{F}} f_\af y^*_\af
\leq \sum_{\af\in \mathcal{F}} f_\af \hat y_\af$
and $\sum_{\af\in \mathcal{F}} f_\af y^*_\af
\leq \sum_{\af\in \mathcal{F}} f_\af \tilde y_\af$.
By linearity, it holds
\[
f^*_\Dt=
\sum_{\af\in \mathcal{F}} f_\af y^*_\af =
\lmd \sum_{\af\in \mathcal{F}} f_\af \hat y_\af +
(1-\lmd)\sum_{\af\in \mathcal{F}} f_\af \tilde y_\af.
\]
Therefore, we must have $\sum_{\af\in \mathcal{F}} f_\af \hat y_\af=f^*_\Dt$ since
$0<\lmd <1$.
On the other hand, by the definition of $\hat y$,
we know $f(x^*)=\sum_{\af\in \mathcal{F}} f_\af \hat y_\af =f^*_\Dt$.
Thus $x^*$ is one point at which the polynomial $f(x)$
attains its lower bound $f^*_\Dt$,
which implies that $x^*$ is a global minimizer
of $f(x^*)$.
\eproof

\smallskip
The algorithm for minimizing $f(x)$ via sparse
SOS relaxation~\reff{sumsos1}-\reff{sumsos3}
is as follows.
\begin{alg}[Minimizing sum of polynomials]
\qquad
\bdes
\item [Input:] $n,m, \Dt_i,f_i(x_{\Dt_i})\,(i=1,\cdots,m)$

\item [Output:] $\mathcal{V}_i$ and $\mathcal{X}_{\Dt_i}\,(i=1,\cdots,m)$

\item [Begin] \quad
\bdes
\item [Step~1:] Solve the dual problem \reff{summom1}-\reff{summom3}.
Get the optimal solution $y^*$.
\item [Step~2:] For each $1\leq k\leq n$, find the set $\mathcal{V}_k$
of points that can be extracted from $M_d^{\dt_k}(y^*)$.
\item [Step~3:] For every $k$ with $|\mathcal{V}_k|>1$,
find the set $\mathcal{X}_{\Dt_i}$ from $M_d^{\Dt_i}(y^*)$  whenever $k\in\Dt_i$.
\edes
\item [End]
\edes
\end{alg}

\bigskip

As an example, let us illustrate how to solve the global optimization problem
\[
\min_{x\in\re^3} \quad  \underbrace{(x_1^2-1)^2 + (x_1-x_2)^4}_{f_1(x_{\Dt_1})} +
\underbrace{(x_2-x_3)^4}_{f_2(x_{\Dt_2})}
\]
and find global minimizers.
Here $\Dt_1=\{1,2\}$ and $\Dt_2=\{2,3\}$.
Solve the dual problem~\reff{summom1}-\reff{summom3}
and we get solutions
\[
\mathcal{M}_1^{\Dt_1}(y^*) = \mathcal{M}_1^{\Dt_2}(y^*)=
\bbm
    1 & 0 & 0 & 1 & 1 & 1 \\
    0 & 1 & 1 & 0 & 0 & 0  \\
    0 & 1 & 1 & 0 & 0 & 0  \\
    1 & 0 & 0 & 1 & 1 & 1  \\
    1 & 0 & 0 & 1 & 1 & 1  \\
    1 & 0 & 0 & 1 & 1 & 1  \\
\ebm.
\]
Both $\mathcal{M}_1^{\Dt_1}(y^*)$ and $\mathcal{M}_1^{\Dt_2}(y^*)$
have rank two and satisfy the flat extension condition.
Using the technique from \cite{HenLas}, we can extract
\[
\mathcal{V}_1 =
\mathcal{V}_2 =
\mathcal{V}_3 = \{-1,1\}
\]
and
\[
\mathcal{X}_{\Dt_1} = \left\{ \bbm -1 \\ -1\ebm, \bbm 1 \\ 1 \ebm  \right\}, \quad
\mathcal{X}_{\Dt_2} = \left\{ \bbm -1 \\ -1\ebm, \bbm 1 \\ 1 \ebm  \right\}.
\]
Since the $x_2$-component from $\mathcal{X}_{\Dt_1}$ and $\mathcal{X}_{\Dt_2}$
must be the same, we know there are two global minimizers
$x^*=\pm (1,1,1)$.

\subsection{Nonlinear least squares problems}

Now we consider the special case that each $f_i(x_{\Dt_i})$
is a square of some polynomial, say, $f_i(x_{\Dt_i})=g_i^2(x_{\Dt_i})$.
Then the global minimization of $f(x)=\sum_i\, f_i(x_{\Dt_i})$
is equivalent to solving the nonlinear least squares (NLS) problem
associated with the polynomial system:
\begin{align} \label{nls}
g_1(x_{\Dt_1})=g_2(x_{\Dt_2})=\cdots=g_m(x_{\Dt_m})=0.
\end{align}
In this situation, the polynomial function
is often nonconvex and it is very difficult
for general numerical optimization schemes like branch-bound
to find the global minimizer of $f(x)$.

\begin{thm}\label{nlssos}
If the polynomial system \reff{nls} admits a solution,
then the sparse SOS relaxation \reff{sumsos1}-\reff{sumsos3}
is exact, i.e., $f^*_\Dt=f^*_{sos}=f^*$.
\end{thm}
\proof
Obviously $f^*=0$.
And $\gamma =0$ is a feasible solution to problem
\reff{sumsos1}-\reff{sumsos3}, since
$f(x)$ itself is a sparse SOS representation
as in \reff{sumsos2}-\reff{sumsos3}.
So $f^*_\Dt\geq 0$, and hence
all the inequalities in the Theorem~\ref{bdrel}
become equalities.
\eproof

\begin{remark}
When the polynomial system \reff{nls} admits
a solution, we necessarily have $f^*=0$.
This might be trivial in some sense.
However, the optimal solution $y^*$ to the dual problem
\reff{summom1}-\reff{summom3}
can help recover the real zeros of
polynomial system \reff{nls},
which are absolutely the global minimizers of $f(x)$.
See the example below.
\end{remark}

\begin{exm}
Consider the sparse polynomial system
\begin{align*}
2x_1^2-3x_1+2x_2-1& = 0 \\
2x_i^2+x_{i-1}-3x_i+2x_{i+1}-1 & = 0 \, (i=2,\cdots,n-1) \\
2x_n^2+x_{n-1}-3x_n-1 & = 0.
\end{align*}
This polynomial system is consistent and has at least two real solutions.
Set $n=20$.
We apply sparse SOS relaxation~\reff{sumsos1}-\reff{sumsos3}
to solve the least squares problem and get the lower bound
$f^*_\Dt \approx -2.0 \cdot 10^{-11}$.
Using the optimal dual solution, we obtain two real solutions
(only the first four digits are shown)
{\scriptsize
\begin{align*}
\hat x &=( 1.8327 ,\,  -0.1097 ,\,  -0.5929 ,\,  -0.6860 ,\, -0.7032 ,\,-0.7064 ,\, -0.7070 ,\,
-0.7071 ,\,  -0.7071 ,\, -0.7071 ,\, \\
& \quad -0.7071 ,\, -0.7070 ,\,-0.7068 ,\,  -0.7064 ,\, -0.7051 ,\,-0.7015 ,\,
 -0.6919 ,\,  -0.6658,\, -0.5960 ,\,  -0.4164) \\
\tilde x & =(-0.5708  ,\, -0.6819 ,\,  -0.7025 ,\,  -0.7063 ,\,  -0.7070 ,\,  -0.7071 ,\,-0.7071 ,\,
-0.7071 ,\,  -0.7071 ,\, -0.7071 ,\, \\
& \quad -0.7071 ,\, -0.7070 ,\, -0.7068  ,\, -0.7064 ,\,-0.7051 ,\,  -0.7015 ,\, -0.6919 ,\,
 -0.6658 ,\, -0.5960  ,\, -0.4164).
\end{align*}
}
\end{exm}

\subsection{A sparser SOS relaxation}

\,From Theorem~\ref{nlssos}, we know
the sparse SOS relaxation~\reff{sumsos1}-\reff{sumsos3}
is exact whenever the polynomial system~\reff{nls} admits
a solution, and the optimal dual solution can help
recover the real zeros.
This fact makes it possible to exploit the sparsity
of each $f_i(x_{\Dt_i})$ further.
In \reff{sumsos1}-\reff{sumsos3}, we assume each
$f_i(x_{\Dt_i})$ is a dense polynomial.
However, if each $f_i(x_{\Dt_i})$ is sparse,
we can get a sparser SOS relaxation.
It is obvious that
\[
\supp(f_i) \subseteq \mathcal{G}_i+\mathcal{G}_i
\]
where $\mathcal{G}_i$ is the convex hull of
$\{ \af\in\N^n:\, 2\af \in \supp(f_i)\}$.
This motivates us to propose the sparser SOS relaxation
\begin{align}
f^*_{\Dt_s} :=\max & \quad \gamma \label{morspar1}\\
s.t. & \quad f(x) - \gamma = \sum_{i=1}^m
\mathbf{m}_{\mathcal{G}_i}(x)^T W_i
\mathbf{m}_{\mathcal{G}_i}(x)  \label{morspar2} \\
& \qquad W_i\succeq 0, \, i=1,\cdots,m. \label{morspar3}
\end{align}
Here $\mathbf{m}_{\mathcal{G}_i}(x_{\Dt_i})$ is the column
vector of all monomials in $x$
with exponents from $\mathcal{G}_i$.
The size of matrix $W_i$ is equal to the cardinality of
$\mathcal{G}_i$.
Similar to \reff{sumsos1}-\reff{sumsos3},
the dual of \reff{morspar1}-\reff{morspar3}
can be derived to be
\begin{align}
f^*_{\Sig_s}:=
\min & \quad  \sum_\af f_\af y_\af \label{morsmom1}\\
s.t. & \quad  \mathcal{M}_{\mathcal{G}_i}(y) \succeq 0,
\, i=1,\cdots,m  \label{morsmom2} \\
& \quad y_{0} =1. \label{morsmom3}
\end{align}
Here the sparse moment matrix $\mathcal{M}_{\mathcal{G}_i}(y)$
is indexed by vectors from $\mathcal{G}_i$ and defined as
\[
\mathcal{M}_{\mathcal{G}_i}(y)(\af,\bt) = y_{\af+\bt}
\]
for all $\af,\bt\in\mathcal{G}_i$.

\begin{thm}\label{morebdrel}
The optimal values $f^*_{\Sig_s},\,f^*_{\Dt_s},\,f^*_{\Sig},\,f^*_\Dt,\,f^*_{sos},\, f^*$
satisfy the relationship
\[
f^*_{\Sig_s} = f^*_{\Dt_s} \leq f^*_\Sig=f^*_\Dt \leq f^*_{sos} \leq f^*.
\]
\end{thm}
\proof
Applying the standard duality theory in convex programming
as in the proof of Theorem~\ref{bdrel},
we can get the first equality from the left
by proving \reff{morspar2}-\reff{morspar3}
has a strict interior point.
Since the relaxation \reff{morspar1}-\reff{morspar3}
is a special case of \reff{sumsos1}-\reff{sumsos3},
we obtain the first inequality from the left.
The other relations follow Theorem~\ref{bdrel}.
\eproof

\begin{thm}\label{sparnls}
Suppose $f_i(x_{\Dt_i})=g_i^2(x_{\Dt_i})$.
If the polynomial system \reff{nls} admits a solution,
then the sparse SOS relaxation \reff{morspar1}-\reff{morspar3}
is exact, i.e., $f^*_{\Dt_s}=f^*_{sos}=f^*$.
\end{thm}
\proof
The proof is almost the same as for Theorem~\ref{nlssos}.
Obviously $f^*=0$.
And $\gamma =0$ is a feasible solution, since
$f(x)$ itself is a sparse SOS representation
as in \reff{morspar2}-\reff{morspar3}.
So $f^*_{\Dt_s}\geq 0$, and hence
all the inequalities in the Theorem~\ref{morebdrel}
become equalities.
\eproof

\begin{remark}
When the polynomial system \reff{nls} admits
a solution, we must have $f^*=0$.
This lower bound itself might not be interesting.
However, the optimal dual solution $y^*$ to
\reff{morsmom1}-\reff{morsmom3}
can help recover the real zeros of
polynomial system \reff{nls},
which are absolutely the global minimizers of $f(x)$.
This observation is very important and has
many applications.
See examples in Subsection~5.1.
\end{remark}

\section{Numerical examples}
\setcounter{equation}{0}

In this section, we present some numerical experiments
using sparse SOS relaxations \reff{sumsos1}-\reff{sumsos3}
and \reff{morspar1}-\reff{morspar3}.
First, we use them to solve some test problems
from unconstrained optimization.
Second, we generate various random polynomials,
test the performance of these sparse SOS relaxations
and compare with other methods.
All the computations are implemented
on a Linux machine with 0.98 GB memory and 1.46 GHz CPU.
The SOS relaxations are solved by softwares
{\it SeDuMi} \cite{sedumi} and {\it YALMIP} \cite{yalmip}.
Throughout this section, the computation time is in CPU seconds.
The accuracy of relaxations is measured by
$\frac{|f(\hat x) - \hat f|}{\max\{1,|f(\hat x)|\}}$,
where $\hat x$ is one extracted solution
and $\hat f$ is the computed lower bound.

\subsection{Some global optimization test problems}

In this subsection, we apply SOS relaxations \reff{sumsos1}-\reff{sumsos3}
and \reff{morspar1}-\reff{morspar3}
to solve some global optimization test problems from
\cite{CGT,MGH,Nash}.
The relaxation \reff{sumsos1}-\reff{sumsos3} is usually applied
when each $f_i(\Dt_i)$ is almost dense, and
the sparser relaxation \reff{morspar1}-\reff{morspar3}
is usually applied when each $f_i(\Dt_i)$ is sparse.
All the test functions in this subsection have
global minimum $f^*=0$.
So we use the absolute value of the lower bounds
$f^*_\Dt$ or $f^*_{\Dt_s}$ to measure the
accuracy of the relaxation.

First, consider the following test functions.
\bit

\item The chained singular function \cite{CGT}:
\[
f(x) = \frac{1}{10^5} \sum_{i\in J} \left(
(x_i+10x_{i+1})^2+5(x_{i+2}-x_{i+3})^2+(x_{i+1}-2x_{i+2})^4+10(x_i-10x_{i+3})^4
\right)
\]
where $J=\{1,3,5,\cdots,n-3\}$ and $n$ is a multiple of $4$.
The factor $\frac{1}{10^5}$ is used to scale the coefficients
to avoid numerical troubles.

\item The chained wood function \cite{CGT}
\begin{align*}
f(x) = \sum_{i\in J} (
100(x_{i+1}-x_i^2)^2+(1-x_i)^2+90(x_{i+3}-x_{i+2}^2)^2+\\
\hspace{1cm} (1-x_{i+2})^2+10(x_{i+1}+x_{i+3}-2)^2+0.1(x_{i+1}-x_{i+3})^2 )
\end{align*}
where $J=\{1,3,5,\cdots,n-3\}$ and $n$ is a multiple of $4$.

\item The generalized Rosenbrock function \cite{Nash}:
\[
f(x) = \sum_{i=2}^n
\left\{ 100(x_i-x_{i-1}^2)^2+(1-x_i)^2 \right\}.
\]

\eit
We apply SOS relaxation \reff{sumsos1}-\reff{sumsos3}
to minimize these polynomial functions.
The relaxation \reff{sumsos1}-\reff{sumsos3} is solved by softwares
{\it SeDuMi} and {\it YALMIP}.
The accuracy
and consumed CPU time are in Table~\ref{LessSpar}.
The problems are solved from size $100$ to $500$.
\begin{table}[htb]
\centering
\begin{tabular}{|c||c|r||c|c||c|c|} \hline
& \multicolumn{2}{c||}{chained singular} & \multicolumn{2}{c||}{chained wood}
& \multicolumn{2}{c|}{gen. Rosen.} \\ \hline
n& accu. & time\, &  accu. & time &  accu. & time \\ \hline
100 & 3.2e-09 & 2.72  &  3.5e-10 &  1.52 & 9.0e-8 & 0.95 \\ \hline
200 & 3.0e-10 & 5.29  &  3.7e-10 &  2.25 & 1.8e-7 & 1.46 \\ \hline
300 & 5.0e-09  & 8.01  & 3.8e-10 &  3.19& 2.7e-7 & 2.24 \\ \hline
400 & 5.0e-10 & 11.64  & 3.9e-10 &  4.12 & 3.6e-7 & 2.88 \\ \hline
500 & 4.9e-09  & 33.09 & 3.9e-10 &  5.12 & 4.5e-7 & 3.45 \\ \hline
\end{tabular}
\caption{The performance of sparse SOS relaxation \reff{sumsos1}-\reff{sumsos3}
}
\label{LessSpar}
\end{table}
For these polynomials,
the relaxation~\reff{sumsos1}-\reff{sumsos3}
is almost the same as the one in \cite{WKKM}.
This is because the csp graphs of these polynomials are chordal graphs.
However, if the csp graphs are sparse but
their chordal extensions are much dense,
then the relaxation in \cite{WKKM}
is very similar to the dense SOS relaxation.
In such situations, the relaxation \reff{sumsos1}-\reff{sumsos3} might be more suitable.
For example, to minimize the sparse polynomial
\[
(x_1^2+x_2^2-1)^2 + (x_2^2+x_3^2-1)^2+\cdots+(x_{n-1}^2+x_n^2-1)^2 +(x_n^2+x_1^2-1)^2,
\]
the chordal extension of the csp graph is the complete graph, and hence
the sparse SOS relaxation using chordal extension
is the same as the dense SOS relaxation.
However, the sparse relaxation \reff{sumsos1}-\reff{sumsos3}
is very suitable for this problem.

\bigskip

Second, consider the following test functions.
\bit
\item Broyden tridiagonal function \cite{MGH}:
\[
f(x) = \sum_{i=1}^n \Big(
(3-2x_i)x_i-x_{i-1}-2x_{i+1}+1 \Big)^2
\]
where $x_0=x_{n+1}=0$.

\item Broyden banded function \cite{MGH}:
\[
f(x) = \sum_{i=1}^n \left(
x_i(2+10x_i^2)+1-\sum_{j\in J_i} (1+x_j)x_j \right)^2
\]
where $J_i=\{j:j\ne i, \max(1,i-5)\leq j\leq \min(n,i+1)\}$.

\item Discrete boundary value function \cite{CGT}:
\[
f(x) = \sum_{i=1}^n \Big(
2x_i-x_{i-1}-x_{i+1}+\half\,h^2(x_i+t_i+1)^3 \Big)^2
\]
where $h=\frac{1}{n+1},\, t_i=ih,\, x_0=x_{n+1}=0$.

\eit
These three polynomials have sparse summand polynomial
$f^*_{\Dt_i}$.
So we apply the sparser SOS relaxation \reff{morspar1}-\reff{morsmom3}
and solve it by using softwares
{\it SeDuMi} and {\it YALMIP}.
The computational results are in Table~\ref{MoreSpar}.
All the problems are solved quite well in a few seconds.

\begin{table}[htb]
\centering
\begin{tabular}{|r|r|r||c|c|r||c|c|r|} \hline
\multicolumn{3}{|c||}{Broyden Tridiagonal} &\multicolumn{3}{c||}{Broyden banded}
&\multicolumn{3}{c|}{disc. bound val. } \\ \hline
 $ n $ & accu. & time &$ n $ & accu. & time &$ n $ & accu. & time \\ \hline
100 & 1.2e-7 & 2.65   & 10  & 3.6e-11 & 9.72 & 10  & 6.0e-12 & 0.92  \\ \hline
200 & 2.3e-7 & 2.69  & 15  & 2.2e-10 & 17.28 &  20  & 3.4e-11 & 1.57  \\ \hline
300 & 5.0e-7 & 3.58  & 20  & 1.6e-10 & 25.27 & 25  & 1.6e-11 & 2.28 \\ \hline
400 & 3.0e-6 & 4.53  & 25 & 1.8e-10 & 35.19 &  30  & 1.1e-11 & 2.47  \\ \hline
500 & 4.1e-6 & 5.44  & 30 & 4.9e-10 & 45.30 &  35  & 3.9e-11 & 3.00  \\ \hline
\end{tabular}
\caption{ The performance of sparse SOS relaxation \reff{morspar1}-\reff{morspar3}
}
\label{MoreSpar}
\end{table}

For the Broyden tridiagonal function,
we can also apply the sparse relaxation~\reff{sumsos1}-\reff{sumsos3}
or chordal extension from \cite{WKKM}.
They are slightly more expensive.
For $n=500$, the problem can be solved in about ten seconds
with similar accuracy.
However, for the Broyden banded function and discrete boundary value function,
the relaxation~\reff{sumsos1}-\reff{sumsos3} and the method in \cite{WKKM}
are much more expensive.
For instance, when $n$ has values $10$ or bigger, they are usually difficult to implement
due to computer memory restrictions.

%

One interesting observation in Table~\ref{MoreSpar} is that
the accuracy for the Broyden tridiagonal function
is not as high as for the other two functions.
One possible reason is that the global minimizer of Broyden tridiagonal function is not unique
and there are additional numerical troubles caused
from extracting minimizers.
This illustrates that
the computation is more numerically difficult
when there are multiple global solutions.

\subsection{Randomly generated test problems}

In this subsection, we present the computational results for
randomly generated polynomials.
The aim is to test the performance of
the sparse SOS relaxation~\reff{sumsos1}-\reff{sumsos3}
for minimizing random polynomials
and compare with other sparse SOS methods.
For these randomly generated polynomials,
solve the sparse relaxation~\reff{sumsos1}-\reff{sumsos3}
by using softwares {\it SeDuMi} and {\it YALMIP}.
Then we get lower bounds $f^*_\Dt$ and
extract minimizers $\hat x$.
Since we do not know the true global minimizers in advance,
the accuracy of $\hat x$ can be measured by
$err\, = \frac{|f(\hat x) - f^*_\Dt|}{\max\{1,|f(\hat x)|\}}.$
The smaller $err$ is, the more accurate $\hat x$ is,
since $f^*_\Dt$ is a guaranteed lower bound.

\bigskip
\subsubsection{Randomly generated sums of small polynomials}
\bigskip

In this subsubsection,
we randomly generate sparse polynomials $f(x)$ of the form~\reff{sumpoly}
and use them to test the performance of
the sparse relaxation~\reff{sumsos1}-\reff{sumsos3}.
Then the csp graph of $f(x)$ is usually not chordal
and its chordal extension is often much less sparse.
So the method in [34] is usually expensive for these polynomials.
We let $m=n$ and choose $f_i$ to have form
\[
f_i(x_{\Dt_i}) = \mathbf{m}_d(x_{\Dt_i})^T \cdot A_i \cdot \mathbf{m}_d(x_{\Dt_i})
+ b_i^T \mathbf{m}_{2d-1}(x_{\Dt_i})
\]
where $\Dt_i$ are chosen to be random subsets of $[n]$
with cardinality at most $\|\Dt\|$.
Here $N_{i}= \binom{|\Dt_i|+d}{d}$,
$A_i = n I_{N_{i}} + BB^T$,
$B \in \re^{N_{i}\times N_{i}}$
and $b_i \in \re^{\binom{ |\Dt_i|+d-1}{d-1} } $
are random.
So each $A_i$ is positive definite.
This choice guarantees that the global minimizers of $f(x)$
are contained in some compact set.

\begin{table}[htb]
\centering
\begin{tabular}{||c||c|c|c||c||c|c|c||c||} \hline
  & \multicolumn{4}{c||}{$\|\Dt\|=3$} &
\multicolumn{4}{c||}{$\|\Dt\|=4$}  \\ \hline
  & \multicolumn{3}{c||}{CPU seconds} & accu & \multicolumn{3}{c||}{CPU seconds} & accu \\ \hline
$n$  & max & avr. & min & max & max & avr. & min & max \\ \hline
20 & 0.85 & 0.62 & 0.54 & 4.1e-9 & 1.46 & 1.15 & 0.91 & 2.4e-9 \\ \hline
40 & 1.22 & 1.07 & 0.91 & 1.9e-9 & 2.86 & 2.49 & 2.25 & 2.9e-9  \\ \hline
60 & 1.80 & 1.55 & 1.45 & 2.9e-9 & 4.43 & 4.17 & 3.91 & 3.1e-9  \\ \hline
80 & 2.30 & 2.18 & 2.02 & 2.3e-9 & 6.26 & 5.94 & 5.24 & 3.7e-9  \\ \hline
100 & 3.02 & 2.70 & 2.33 & 2.8e-9  & 7.85 & 7.41 & 7.01 & 5.0e-9 \\ \hline
\end{tabular}
\caption{Computational results for quartic polynomials with different sizes}
\label{degfour}
\end{table}

First, let $2d=4$ and $n$ be $20,40,60,80,100$.
For each $\|\Dt\|$ (3 or 4) and $n$,
we generate $100$ random polynomials in the way described above.
For each one, we solve the sparse SOS relaxation \reff{sumsos1}-\reff{sumsos3}
by using softwares {\it SeDuMi} and {\it YALMIP},
and get the lower bound $f^*_\Dt$ and optimal dual solution $\hat y$.
For all these randomly generated polynomials,
the moment matrices $\mathcal{M}_d^{\Dt_i}(\hat y)$ have numerical rank one.
So we can easily extract the minimizer $\hat x$.
The maximum, average and minimum of consumed CPU time
are in Table~\ref{degfour}.
For these random polynomials, we just record the
maximum error of the extracted minimizers. From Table~\ref{degfour},
for $\Dt=4$,
we can find the global minimizer of a quartic sparse polynomial of $100$
variables with error $O(10^{-9})$ within about $8$ CPU seconds.

\begin{table}[htb]
\centering
\begin{tabular}{||c||c|c|c||c||c|c|c||c||} \hline
  & \multicolumn{4}{c||}{$\|\Dt\|=3$} &
\multicolumn{4}{c||}{$\|\Dt\|=4$}  \\ \hline
  & \multicolumn{3}{c||}{CPU seconds} & err & \multicolumn{3}{c||}{CPU seconds} & err \\ \hline
$2d$  & max & avr & min & max & max & avr & min & max \\ \hline
4 & 1.01 & 0.87 & 0.77 & 2.5e-9 & 2.33 & 1.93 & 1.65  & 2.4e-9 \\ \hline
6 & 3.22 & 2.96 & 2.67 & 1.8e-9 & 17.16  & 14.92 & 11.71 & 2.2e-9 \\ \hline
8 & 13.07 & 11.44 & 10.13 & 1.7e-8 & 136.67  & 119.90 & 107.28 & 9.4e-8 \\ \hline
\end{tabular}
\caption{Computational results for polynomials of size $n=30$
with different degrees}
\label{size30}
\end{table}

Second, let $n=30$ and $2d$ be $4,6,8$.
For each $\|\Dt\|$ (3 or 4) and $2d$,
we generate $100$ random polynomials in the way described in the above.
For each one, solve the sparse SOS relaxation \reff{sumsos1}-\reff{sumsos3}
by using softwares {\it SeDuMi} and {\it YALMIP},
and get the lower bounds $f^*_\Dt$ and optimal dual solution $\hat y$.
Similarly, all moment matrices $\mathcal{M}_d^{\Dt_i}(\hat y)$ have rank one,
and the minimizer $\hat x$ can be extracted easily.
The maximum, average and minimum of the consumed CPU time, and the maximum
error of extracted minimizers are in Table~\ref{size30}. For $\|\Dt\|=4$,
the global minimizer of such generated polynomials of degree $8$
and $30$ variables can be found with error $O(10^{-8})$ within about $120$ seconds.

\medskip

We remark that the sparsity technique in \cite{WKKM}
is too expensive to be implementable
for minimizing these random polynomials generated in the way as above
because of computer memory limitations.
For these random polynomials,
the sparsity technique using chordal extension
is almost as expensive as the general dense SOS relaxation.
This is because the chordal extensions of csp graphs of these polynomials
are usually much more dense than the original csp graphs.
However, as we have seen in the above,
the SOS relaxation~\reff{sumsos1}-\reff{sumsos3}
is very suitable for these polynomials.

%

\bigskip
\subsubsection{Random sparse polynomials with given chordal extension}
\bigskip

In this subsubsection, we generate random sparse polynomials
in a similar way as in \cite{WKKM}, and compare the performance
of our sparse SOS relaxation \reff{sumsos1}-\reff{sumsos3}
with the one in \cite{WKKM} using chordal extension.
Generate a chordal graph randomly such that the size
of every maximal clique is at most $6$.
Let $\{C_1,\cdots,C_m\}$ be the set of maximal cliques.
If we choose $\Dt_i = C_i$,
then the sparse SOS relaxation \reff{sumsos1}-\reff{sumsos3}
is the same as the one using chordal extension.
Therefore, to make a reasonable comparison,
for each $C_i$, we choose a random subset $\Dt_i\subseteq C_i$.
Choose each small polynomial $f_i$ to have the form
\[
f_i(x_{\Dt_i}) = \mathbf{m}_d(x_{\Dt_i})^T \cdot A_i \cdot \mathbf{m}_d(x_{\Dt_i})
+ b_i^T \mathbf{m}_{2d-1}(x_{\Dt_i}).
\]
Here $N_{i}= \binom{|\Dt_i|+d}{d}$,
$A_i = n I_{N_{i}} + BB^T$,
$B \in \re^{N_{i}\times N_{i}}$
and $b_i \in \re^{\binom{ |\Dt_i|+d-1}{d-1} } $
are random.
The global minimizers of $f(x)=\sum_i f_i(x_{\Dt_i})$ generated as above always exist and
are contained in some compact set.

\begin{table}[htb]
\centering
\begin{tabular}{|c|c|c|c|c|c|c|c|c|} \hline
  & \multicolumn{4}{|c|}{ relaxation \reff{sumsos1}-\reff{sumsos3} } &
\multicolumn{4}{|c|}{ relaxation using chordal extension }  \\ \hline
  & \multicolumn{3}{c|}{CPU seconds} & accu & \multicolumn{3}{c|}{CPU seconds} & accu \\ \hline
$n$  & max & avr. & min & max & max & avr. & min & max \\ \hline
20 & 1.75  & 1.21  & 0.96  & 6.8e-9 & 2.15  & 1.78  & 1.43  & 5.5e-9 \\ \hline
40 & 3.07  & 2.69  & 2.24  & 7.5e-9 & 4.08  & 3.51 & 3.12   & 4.9e-9  \\ \hline
60 & 4.99  & 4.54  & 3.82  & 6.7e-9 & 7.88  & 6.93  & 5.65  & 6.4e-9  \\ \hline
80 & 6.59  & 5.87  & 5.23  & 6.3e-9 & 10.84  & 9.57  & 8.59   & 5.7e-9  \\ \hline
100& 9.34  & 7.64  & 7.11  & 7.2e-9  & 13.45  & 12.76  & 11.74  & 4.3e-9 \\ \hline
\end{tabular}
\caption{Comparison with chordal extension
on quartic polynomials }
 \label{chordextdeg4}
\end{table}

For polynomials randomly generated as above,
the technique in \cite{WKKM} using chordal extension is a good choice,
because there exists one sparse chordal extension of the csp graph.
Now we compare the computational results
for these two methods.

\begin{table}[htb]
\centering
\begin{tabular}{|c|c|c|c|c|c|c|c|c|} \hline
  & \multicolumn{4}{|c|}{ relaxation \reff{sumsos1}-\reff{sumsos3} } &
\multicolumn{4}{|c|}{ relaxation using chordal extension }  \\ \hline
  & \multicolumn{3}{c|}{CPU seconds} & accu & \multicolumn{3}{c|}{CPU seconds} & accu \\ \hline
$2d$  & max & avr. & min & max & max & avr. & min & max \\ \hline
4 & 2.87  & 1.98  & 1.35  & 7.2e-9 & 3.06  & 2.21  & 1.69   & 4.3e-9  \\ \hline
6 & 22.61   & 16.78  & 10.53  & 6.9e-9 & 32.15   & 23.91   & 13.51  & 5.1e-9  \\ \hline
8 & 193.45  & 131.17  & 98.75  & 6.7e-9 & 253.79  & 186.84  & 112.37  & 5.8e-9  \\ \hline
\end{tabular}
\caption{Comparison with chordal extension
on polynomials with $30$ variables }
\label{chordextn30}
\end{table}

First, let $2d=4$ and $n$ be $20,40,60,80,100$.
For each $n$,
generate $50$ random polynomials as above.
For each of these random polynomials, 
solve the relaxation \reff{sumsos1}-\reff{sumsos3},
find a chordal extension of the csp graph of $f(x)$
and then apply the sparse relaxation in \cite{WKKM}.
Both relaxations are solved
by softwares {\it SeDuMi} and {\it YALMIP}.
Then we extract minimizers $\hat x$ from moment matrices.
The computational results
are in Table~\ref{chordextdeg4}.
For these solved problems, we just record the
maximum error of the relaxation.
Second, let $n=30$ and $2d$ be $4,6,8$.
For each $2d$,
generate $50$ random polynomials as above.
For each one, solve the problem by
the relaxation \reff{sumsos1}-\reff{summom3}
and the one in \cite{WKKM}
using chordal extension.
They are solved by softwares {\it SeDuMi} and {\it YALMIP}.
The computational results are in Table~\ref{chordextn30}.

\,{From} Tables \ref{chordextdeg4} and \ref{chordextn30},
we observe that for polynomials randomly generated as above
the sparse SOS relaxation~\reff{sumsos1}-\reff{sumsos3}
is slightly more computationally efficient
than the one using chordal extension.
As we can see, for these random polynomials,
there is not much difference between the qualities of these two kinds of sparse SOS relaxations.
The distinction between their qualities depends on specific problems.
Of course, theoretically the sparse relaxation
using chordal extension in \cite{WKKM} is at least as tight as
the relaxation~\reff{sumsos1}-\reff{sumsos3}.

\bigskip
\subsubsection{Random dense polynomials}
\bigskip

In this subsubsection, we test the performance of
our sparse SOS relaxation
on minimizing general dense polynomials.
We observe that every polynomial $f(x)$ is a summation of monomials
whose number of variables is at most the degree $\deg(f)$.
So the sparse SOS relaxation \reff{sumsos1}-\reff{sumsos3}
is attractive when the degree $2d$ is small like $4$.
We generate the random dense polynomials as follows
\[
f(x) = \mathbf{m}_d(x)^T \cdot A \cdot \mathbf{m}_d(x)
+ b^T \mathbf{m}_{2d-1}(x).
\]
Here $N = \binom{n+d}{d}$,
$ A = n I_{N} \, + BB^T$,
$B\in \re^{N\times N}$ is a random matrix
and $b\in \re^{\binom{n+2d-1}{n}}$ is a random vector.
So the global minimizers of $f(x)$
generated this way
are contained in some compact set.
Note that $f(x)$ is also a summation of small polynomials.
Let $\Dt_i$ be the subsets of $[n]$ with cardinality $2d$.
Then we can write
\[
f(x) = \sum_{i=1}^{\binom{n}{2d}} f_i(x_{\Dt_i})
\]
for some small polynomials $f_i(x_{\Dt_i})$.

%
%
%

\begin{table}[htb]
\centering  \footnotesize
\begin{tabular}{|r|c|c|c|c|c|c|c|c|c|} \hline
      $n$ & 16 & 17 & 18 & 19 & 20 & 21 & 22  & 23 \\  \hline
max  time &  335.29   & 569.74   & 901.32   &  1505.45  & 2249.19   & 3257.86   & 4734.25 & 7060.72 \\  \hline
avr. time &  241.48   & 455.32   & 751.69   &  1245.22  & 2070.70   & 2989.45   & 4497.84 & 6419.53  \\  \hline
min  time &  205.60   & 397.11   & 688.58   &  1052.70  & 1893.02   & 2676.62   & 4197.95 & 5874.28 \\  \hline
accuracy  &  7.3e-9   & 6.7e-9   & 7.4e-9
 & 6.9e-9    & 8.1e-9   & 6.5e-9   &  7.9e-9  & 8.5e-9 \\  \hline
\end{tabular}
\caption{Computational results for dense quartic polynomials}
 \label{densepoly46}
\end{table}

Since $\|\Dt\|=2d$, which should not be big for the effectiveness
of the sparse relaxation \reff{sumsos1}-\reff{sumsos3},
we test for the case that $2d=4$.
Let $n$ be $16,17,18,19,20,21,22,23$.
For each pair $(n,d)$ of these values,
generate $50$ random examples as above.
For each random polynomial,
solve the sparse relaxation \reff{sumsos1}-\reff{sumsos3}
by using softwares {\it SeDuMi} and {\it YALMIP}.
The consumed CPU time
and the accuracy of relaxation are
in Table~\ref{densepoly46}.
We can see that the obtained solutions
are very good within reasonably acceptable time.
When $n\geq 24$, the sparse relaxation~\reff{sumsos1}-\reff{sumsos3}
is then also too expensive to be implementable due to
computer memory restrictions.

\smallskip

For these randomly generated dense polynomials,
the general dense SOS relaxation
and sparse SOS relaxations like in \cite{WKKM}
are not implementable for $n\geq 16$,
due to either computer memory shortage
or the consumed time being more than $10$ hours.
However, when $2d$ is small like $4$,
the sparse SOS relaxation \reff{sumsos1}-\reff{sumsos3}
can solve bigger dense polynomial optimization problems
which can not be solvable by other methods.

%
%
%
%
%
%
%
%
%
%
%
%
%
%
%
%
%
%
%
%
%
%
%
%
%
%
%
%
%
%
%
%
%
%
%
%
%
%
%
%
%
%
%
%
%
%
%
%
%
%
%

\section{Applications}
\setcounter{equation}{0}

Minimizing a summation of small polynomials arises in various applications.
Many big polynomials in applications often come in this form.
In such situations, the sparse SOS relaxation
\reff{sumsos1}-\reff{sumsos3} or \reff{morspar1}-\reff{morspar3}
is very useful.
In this section, we show some applications
in solving sparse polynomial systems
and sensor network localization.

\subsection{Solving sparse polynomial system}

Suppose we are trying to solve the sparse polynomial system
\begin{align*}
g_1(x_{\Dt_1}) =0,\,
 g_2(x_{\Dt_2}) =0,\, \cdots, \,
g_m(x_{\Dt_m}) =0.
\end{align*}
In some applications, these equations are redundant
or even inconsistent.
When the polynomial system does not admit a solution,
we want to seek a least squares solution,
which is often useful in applications.

This problem can be formulated as
finding the global minimizer of the sparse polynomial
\[
f^*:=\min_{x\in\re^n}\, f(x)=\sum_{i=1}^m\, g_i^2(x_{\Dt_i}).
\]
The polynomial system has a real zero if and only if $f^*=0$.
When $f^*=0$, the global minimizers are precisely the
real zeros of the polynomial system.
When $f^*>0$, the global minimizers are the least squares
solutions.

%
%

One important sparse polynomial system of the above form
is from computing the numerical solutions of nonlinear differential
equations. Consider the two-point boundary value problem (BVP)
\[
F(t,x,x^{\prm},x^{\prm\prm})=0,\,\,\, x(a)=\af, \, x(b)=\bt
\]
where $F(t,x,x^\prm,x^{\prm\prm})$ is polynomial function
in $t,x,x^\prm,x^{\prm\prm}$.
To find the numerical solution, the central difference approximation
with a uniform mesh is often used to discretize the derivatives.
Let $N$ be a positive integer and set $h=\frac{b-a}{N+1}$.
Then we get polynomial difference equations
\[
F\left(t_k,x_k, \frac{x_{k+1}-x_{k-1}}{2h},
\frac{x_{k-1}-2x_k+x_{k+1}}{h^2}
\right)=0, \, k=1,\cdots,N\\
\]
where $x_0=\af,\,x_{N+1}=\bt$ and $t_k=a+hk$.
Every polynomial on the left involves $2$ or $3$ variables
$x_{k-1},\,x_k,\,x_{k+1}$.
So this is a sparse polynomial system.
There are several methods for solving this kind of polynomial
system, like Newton's method and homotopy methods.
Newton's method is very fast, but often require an accurate initial guess.
Homotopy methods do not require a ``satisfactory'' guess and work
well for small $N$, but are expensive to implement for large $N$.
We refer to \cite{ABSW} and the references therein
for work in this area.
When $N$ is large, this polynomial system is large but sparse.
We solve this system by
applying the sparse SOS relaxation~\reff{morspar1}-\reff{morsmom3}
for big $N$ (up to $100$ or even bigger).

\begin{exm}[\cite{ABSW}]\label{bvp1}
Consider a basic BVP
\[
x^{\prm\prm} - 2x^3=0,\,\,\, x(0)=\half, \, x(1)=\frac{1}{3}.
\]
The exact solution to this problem is $x(t)=\frac{1}{t+2}$.
Now we discretize the differential equation with mesh size
$h=\frac{1}{N+1}$, then get the difference equation
\begin{align*}
\half-2x_1+x_2 - 2h^2 x_1^3 & = 0, \\
x_{k-1}-2x_k+x_{k+1} - 2h^2 x_k^3 & = 0, \, k=2,\cdots,N-1\\
x_{N-1}-2x_N+\frac{1}{3} - 2h^2 x_N^3 & = 0.
\end{align*}
This is a polynomial system about $x_1,x_2,\cdots,x_N$.
We can solve this polynomial system as a nonlinear least squares problem
by applying sparse SOS relaxation~\reff{morspar1}-\reff{morsmom3}.
The computational results are in Table~\ref{odexamp1}.
\begin{table}[htb] \small
\centering
\begin{tabular}{|r|c|c|c|r|} \hline
   N\, & eqn. error & $\| x_k -x(t_k)\|_\infty$ & $\| x_k -x(t_k)\|_\infty /h^2$  & time   \\ \hline
   5 & 2.8937e-07 & 7.0252e-05 &  2.5291e-003  & 0.52  \\ \hline
   10 & 2.3329e-07 & 1.9570e-05 & 2.3680e-003  & 0.77 \\ \hline
   20 & 5.2879e-07 & 1.5041e-05 & 6.6331e-003  & 1.18 \\ \hline
   30 & 2.6194e-07 & 1.9413e-05 & 1.8656e-002  & 2.09 \\ \hline
   40 & 3.0304e-07 & 4.3344e-05 & 7.2861e-002  & 3.99 \\ \hline
   50 & 6.5375e-07 & 1.5124e-04 & 3.9338e-001  & 6.82 \\ \hline
   60 & 1.5271e-06 & 4.8695e-04 & 1.8119e+00  & 7.77  \\ \hline
   70 & 1.2555e-06 & 5.2428e-04 & 2.6429e+00  & 9.16  \\ \hline
   80 & 9.7315e-07 & 6.1330e-04 & 4.0239e+00  & 9.78  \\ \hline
   90 & 2.7519e-06 & 1.9311e-03 & 1.5991e+01  & 10.81 \\ \hline
   100 & 1.8628e-06 & 8.1425e-04 & 8.3062e+00  & 8.79 \\ \hline
\end{tabular}
\caption{The performance of \reff{sumsos1}-\reff{summom3}
solving the equations in Example~\ref{bvp1}}
\label{odexamp1}
\end{table}
The equation error is defined to be the infinity norm of
the residuals of the left hand side of the polynomial system,
which measures the quality of how the polynomial systems are solved.
The obtained solutions have equation error from $\mathcal{O}(10^{-6})$ to
$\mathcal{O}(10^{-7})$.
If we want to make them more accurate,
they can be used as the initial guesses in Newton's methods for refining.
The accuracy of the discretization is defined to be
the difference between computed solution $x_k$ and
true solution $x(t_k)=\frac{1}{2+t_k}$ where $t_k=\frac{k}{N+1}$.
Since the discretization has error $\mathcal{O}(h^2)$,
we expect that $\| x_k -x(t_k)\|_\infty /h^2$ is a constant.
When $N\leq 40$, we can see that
$\| x_k -x(t_k)\|_\infty /h^2$ is almost constant.
When $N\geq 50$, {\it SeDuMi} experienced numerical troubles,
and the returned solutions are not as accurate as for the smaller $Ns$.
This explains why $\| x_k -x(t_k)\|_\infty$ and $\| x_k -x(t_k)\|_\infty /h^2$
becomes bigger for $N\geq 50$.
Time records the CPU seconds consumed by the SDP solver {\it SeDuMi}.
For $N=100$, the computation takes less CPU time than
for $N=80$ or $N=90$.
This is because the numerical troubles make {\it SeDuMi} terminate
earlier.
\end{exm}

\begin{exm}
Consider another BVP
\[
x^{\prm\prm} + \half (x+t)^3=0,\,\,\, x(0)=0, \, x(1)=0.
\]
Now we discretize the differential equation with mesh size
$h=\frac{1}{N+1}$, then get the difference equation
\begin{align*}
2x_1-x_2 +\half h^2 (x_1+t_1)^3 & = 0  \\
2x_i-x_{i-1}-x_{i+1} +\half h^2 (x_i+t_i)^3 & =0, \, i=2,\cdots,N-1\\
2x_N-x_{N-1}+\half h^2 (x_N+t_N)^3 &=0
\end{align*}
This is a polynomial system about $x_1,x_2,\cdots,x_N$.
We can solve this polynomial system as a nonlinear least squares problem
by applying sparse SOS relaxation~\reff{morspar1}-\reff{morsmom3}.
When $N=30$, we get the following real solution within about $2.5$ CPU seconds
(only the first four digits are shown):
{\scriptsize
\begin{align*}
(-0.0159 ,\,  -0.0312  ,\, -0.0459  ,\, -0.0600  ,\, -0.0735 ,\,
  -0.0864 ,\,  -0.0985 ,\, -0.1099  ,\, -0.1205 ,\,  -0.1302 ,\, \\
-0.1391 ,\,  -0.1470 ,\,  -0.1540 ,\,  -0.1599  ,\,  -0.1646 ,\,
  -0.1682 ,\,  -0.1705  ,\, -0.1715 ,\,  -0.1710 ,\,  -0.1689 ,\, \\
\quad -0.1651 ,\,   -0.1596 ,\,  -0.1521 ,\,  -0.1425 ,\,  -0.1307 ,\,
  -0.1164 ,\,  -0.0995  ,\, -0.0796 ,\, -0.0567 ,\,  -0.0302 ).
\end{align*}
}
\end{exm}

\subsection{Sensor Network Localization}

The {\it sensor network location} problem is basically described as follows:
find a sequence of unknown vectors
$x_1,x_2,\cdots,x_n \in \re^k(k=1,2,\cdots)$
(they are called {\it sensors})
such that distances between these sensors and some other known vectors
$a_1,\cdots,a_m$ (they are called {\it anchors})
are equal to some given numbers.
Now each $x_i$ itself is a $k$-dimensional vector.
To be more specific,
let $\mathcal{A}=\{(i,j)\in [n]\times [n]:\, i<j,\, \|x_i-x_j\|_2 = d_{ij}\}$,
and $\mathcal{B}=\{(i,k)\in [n]\times [m]:\,\, \|x_i-a_k\|_2 = e_{ik}\}$,
where $d_{ij},e_{ik}$ are given distances.
Then the sensor network localization problem
is to find vectors $x_1,x_2,\cdots,x_n$
such that $\|x_i-x_j\|_2 = d_{ij}$ for all $(i,j)\in \mathcal{A}$
and $\|x_i-a_k\|_2 = e_{ik}$ for all $(i,k)\in \mathcal{B}$.
Notice that $\mathcal{A}$ and $\mathcal{B}$ only give some
partial pairs of distances. $\mathcal{A}$ does not contain all the pairs
$(i,j)$ such that $i<j$, and neither does $\mathcal{B}$.

Sensor network localization is also known as
the {\it graph realization problem} or
the {\it distance geometry problem}.
Given a graph $G=(V,E)$ along with a real number associated with each edge,
graph realization is to assign each vertex a coordinate so that
the Euclidean distance between any two adjacent vertices
is equal to the real number associated with that edge.

The locations of sensors can be determined from the polynomial system
\begin{align*}
\| x_i - x_j \|_2^2 &= d_{ij}^2, \forall\, (i,j) \in\mathcal{A},  \\
\| x_i - a_k \|_2^2 &= e_{ik}^2, \forall\, (i,k) \in\mathcal{B} .
\end{align*}
Usually solving this polynomial system directly is very expensive.
Here we solve this polynomial system as a nonlinear least squares problem.
Minimize the quartic polynomial function
\begin{align} \label{L2err}
f(x):=
\sum_{(i,j)\in \mathcal{A}}
\left( \| x_i - x_j \|_2^2 - d_{ij}^2\right)^2 +\,
\sum_{(i,k)\in \mathcal{B}}
\left( \| x_i - a_k \|_2^2 - e_{ik}^2 \right)^2.
\end{align}
where $x=[\,x_1,\,\cdots,\,x_n]$.
$x^*$ is a solution to sensor network localization problem
if and only if $x^*$ is a global minimizer of $f(x)$ such that
$f(x^*)=0$. When $x^*$ is a global minimizer such that $f(x^*) > 0$,
the distances $d_{i,j}$ and $e_{ik}$ are not consistent,
and $x^*$ is a solution in the least squares sense.
This polynomial $f(x)$ is of the form~\reff{sumpoly},
and our sparse SOS relaxation~\reff{sumsos1}-\reff{sumsos3}
can be applied to solve the problem.

We randomly generate test problems
which are similar to those given in \cite{BLTWY05}.
First, we randomly generate $n=500$ sensor locations
$x_1^*,\cdots,x_n^*$ from the unit square
$[-0.5,\,0.5]\times [-0.5,\,0.5]$.
The anchors $\{a_1,a_2,a_3,a_4\}\, (m=4)$ are chosen to be
four fixed points $(\pm 0.45,\, \pm 0.45)$.
Choose edge set $\mathcal{A}$ such that
for every sensor $x_i^*$ there are at most $10$ sensors
$x_j^*\,(j>i)$ with $(i,j)\in \mathcal{A}$ and
$\|x_i^*-x_j^*\|_2 \leq 0.3$.
For every $(i,j)\in \mathcal{A}$,
compute the distance $\|x_i^*-x_j^*\|_2 = d_{ij}$.
Choose edge $\mathcal{B}$ such that every anchor
is connected to all the sensors within distance $0.3$.
For every $(i,k)\in \mathcal{B}$,
compute the distance $\|x_i^*-a_k\|_2 = e_{ik}$.
Then we apply sparse SOS relaxation \reff{sumsos1}-\reff{summom3}
to minimize polynomial function~\reff{L2err}.
The accuracy of computed sensor locations
$\hat x_1,\cdots, \hat x_n$
will be measured by the Root Mean Square Distance (RMSD)
which is defined as
$
\mbox{RMSD} = \left(\frac{1}{n}
\sum_{i=1}^n\|\hat x_i -x_i^* \|_2^2 \right)^{\frac{1}{2}}.
$
We use {\it SeDuMi} to solve the sparse SOS
relaxation~\reff{sumsos1}-\reff{summom3}
on a Linux machine with 1.46 GHz CPU and 0.98GB memory.
The problem can be solved within about $18$ CPU minutes
with accuracy $\mathcal{O}(10^{-6})$.

We refer to \cite{njw_senloc} for more details
about sparse SOS methods for sensor network localization.

\section{Conclusions and discussions}

This paper proposes sparse SOS relaxations for minimizing
polynomial functions that are summations of small polynomials.
We discuss various properties of these relaxations
and the computational issues.
We also present applications of this sparsity technique
in solving polynomial equations
derived from nonlinear differential equations
and sensor network localization.
As a special case,
this sparsity technique provides a heuristic approach to solve
bigger dense polynomial optimization problems.

%

In order to exploit the sparsity,
the polynomial and its SOS representation must be sparse.
In many applications, the polynomials are often
given with sparsity pattern~\reff{sumpoly},
and then the sparsity technique proposed
in this paper is very suitable.
If the sparsity pattern is not given,
one important future work is
how to represent the polynomial in a sparse pattern such that
the technique proposed in this paper is most efficient.
Of course, one simple choice is to consider each monomial
as a small polynomial.

 
The idea of this sparse SOS relaxation
can be applied in a similar way to solve
constrained polynomial optimization problems,
provided the objective and constraint polynomials
are also sums of small polynomials.
See Kim et al. \cite{KKW} and Lasserre~\cite{Las_sparse} for related work.
To get the global minimum, high order
relaxations are usually necessary.
Lasserre~\cite{Las_sparse} proved the convergence
under the running intersection property.
However, unlike the general dense SOS relaxation
for minimizing polynomials over compact sets,
the convergence might fail when
the running intersection property does not hold.
As a counterexample, consider
the {\it Minimum Cover Set Problem}.
Let $G=(V,E)$ be a graph with vertex set $V=[3]$ and
edge set $E=\{(1,2),(1,3),(2,3)\}$.
To find the minimum cover set is equivalent to solving
\begin{align*}
\min_{x\in \re^3} & \quad  f_1(x_{\Dt_1})+f_2(x_{\Dt_2})+f_3(x_{\Dt_3}) \\
s.t. & \quad  x_1^2=x_1,\, x_2^2=x_2, \,x_3^2=x_3, \\
& \quad  x_1+x_2 \geq 1,\, x_1+x_3\geq 1, x_2+x_3\geq 1
\end{align*}
where $\Dt_1=\{1,2\},\,\Dt_2=\{1,3\},\,\Dt_3=\{2,3\}$
and $f_1(x_{\Dt_1})=\half (x_1+x_2),\,
f_2(x_{\Dt_2})=\half (x_1+1_3), f_3(x_{\Dt_3})=\half (x_2+x_3)$.
The running intersection property now fails.
However, we can prove that the global minimum $f^*=2$
and the lower bounds given by sparse SOS relaxations
are at most $\frac{3}{2}$.
The sparse SOS relaxations do not converge for this example.

Another important future work is to apply the sparse SOS relaxations
in solving big real sparse polynomial systems
arising from nonlinear differential equations.

\bigskip

\noindent
{\bf Acknowledge:} The authors wish to thank Bernd Sturmfels and the referees
for helpful suggestions to improve this paper.

\end{document}